\documentclass[3p, final]{elsarticle}
\usepackage{afterpage}

\usepackage{amsmath}
\usepackage{amssymb}
\usepackage{amstext}
\usepackage{amsfonts}

\begin{document}

\begin{frontmatter}

\title{High-Order Resonant Orbit Manifold Expansions For Mission Design In the Planar Circular Restricted 3-Body Problem}

\author[gtaddress]{Bhanu Kumar\corref{mycorrespondingauthor}}
\ead{bkumar30@gatech.edu}
\author[jpladdress]{Rodney L. Anderson}
\cortext[mycorrespondingauthor]{Corresponding author}
\ead{rodney.l.anderson@jpl.nasa.gov}
\author[gtaddress]{Rafael de la Llave}
\ead{rafael.delallave@math.gatech.edu}

\address[gtaddress]{School of Mathematics, Georgia Institute of Technology, 686 Cherry St. Atlanta, GA 30332, USA}
\address[jpladdress]{Jet Propulsion Laboratory/California Institute of Technology, 4800 Oak Grove Dr, Pasadena, CA, 91109, USA}

\begin{abstract}
In recent years, stable and unstable manifolds of invariant objects (such as libration points and periodic orbits) have been increasingly recognized as an efficient tool for designing transfer trajectories in space missions. However, most methods currently used in mission design rely on using eigenvectors of the linearized dynamics as local approximations of the manifolds. Since such approximations are not accurate except very close to the base invariant object, this requires large amounts of numerical integration to globalize the manifolds and locate intersections. 
In this paper, we study hyperbolic resonant periodic orbits in the planar circular restricted 3-body problem, and transfer trajectories between them, by: 1) determining where to search for resonant periodic orbits; 2) developing and implementing a parameterization method for accurate computation of their invariant manifolds as Taylor series; and 3) developing a procedure to compute intersections of the computed stable and unstable manifolds. We develop and implement algorithms that accomplish these three goals, and demonstrate their application to the problem of transferring between resonances in the Jupiter-Europa system.
\end{abstract}


\end{frontmatter}

\section{Introduction}

In recent years, resonant periodic orbits and their stable and unstable manifolds have seen significant interest and use as a tool for trajectory design in multi-body systems. For instance, Anderson and Lo \cite{Anderson2010} demonstrated that a planar version of a Europa Orbiter trajectory designed in 1999 at JPL closely followed stable and unstable manifolds of unstable resonant periodic orbits during resonance transition. They also demonstrated \cite{Anderson2011} the development of new trajectories using homoclinic and heteroclinic connections between resonances. Resonant orbit manifold arcs were also used by Vaquero and Howell \cite{vaqueroHowell} to design transfers from LEO to Earth-Moon libration point orbits. More recently, out of the nine Titan-to-Titan encounters made by Cassini between July 2013 and June 2014, eight of the nine resulting transfers involved resonances \cite{vaqueroCassini}. And even more recently, the baseline mission design for the Europa Lander mission concept made profitable use of these mechanisms for the final approach to the surface of Europa \cite{Anderson2019}. For many other examples of applications of resonant orbits, see Anderson, Campagnola, and Lantoine \cite{Anderson2016}. 

However, the methods used in the previously mentioned studies, as well as in others, rely on using eigenvectors of the linearized dynamics as local approximations of the manifolds. Since such approximations are not accurate except very close to the base invariant object, this requires large amounts of numerical integration to globalize the manifolds and locate intersections, which can decrease accuracy as integration errors add up over longer integration times.

In this paper, we study hyperbolic resonant periodic orbits in the planar circular restricted 3-body problem, and develop methods for accurately computing their manifolds and transfer trajectories between them. We first use the standard Melnikov method \cite{guckholmes} to find Keplerian periodic orbits which survive for small values of the mass parameter $\mu$. This perturbative analysis is followed by numerical continuation to compute the orbits for physically relevant $\mu$ values. We then implement the parameterization method \cite{CabreFontichLlave, haroetal} to compute high order polynomial approximations of the stable and unstable manifolds. Finally, we develop an efficient method which combines the previously computed polynomials with a Poincar\'e section and bisection to compute heteroclinic connections. We also demonstrate application of these tools to the problem of transferring between resonances in the Jupiter-Europa system.

\subsection{Model} \label{modelsection}
The dynamical model considered in the analysis to follow is the well-known planar circular restricted 3-body problem (PCRTBP). In the PCRTBP, one considers two large bodies called the primary body of mass $m_{1}$ and a secondary body of mass $m_{2}$ (collectively referred to as the primaries), revolving about their common center of mass in a circular Keplerian orbit. Units are also normalized so that the distance between the two primaries becomes 1, $\mathcal{G}(m_{1}+m_{2})$ becomes 1, and their period of revolution becomes $2 \pi$. We define a mass ratio $\mu = \frac{m_{2}}{m_{1}+ m_{2}}$, and unless otherwise specified, use a synodic, rotating non-inertial cartesian coordinate system centered at the barycenter of the primaries such that the two primaries are always on the $x$-axis. Due to the normalized units, the primary body will be at $x = -\mu$, and the secondary will be at $x = 1-\mu$. 

\begin{figure}
\centering
\includegraphics[width=0.5\textwidth]{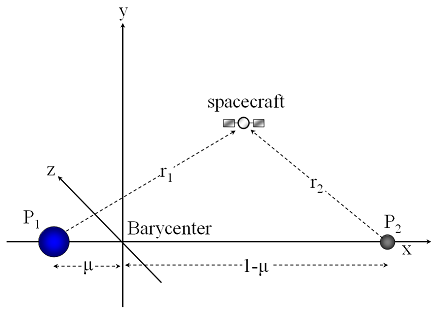}
\caption{\label{fig:crtbp}Diagram of Circular Restricted 3-Body Problem in Synodic Coordinate Frame \cite{ccar}} 
\end{figure}

One then considers the motion of a spacecraft of negligible mass under the gravitational influence of the two primaries. In the planar case we are studying here, we also assume that the spacecraft moves in the same Keplerian orbit plane as the primaries. In this case, and in this synodic coordinate system, the equations of motion become \cite{celletti}
\begin{equation} \label{pcrtbpx} \ddot{x} -2\dot{y}=x -(1-\mu)\frac{x+\mu}{r_{1}^{3}}-\mu \frac{x-1+\mu}{r_{2}^{3}} \end{equation}
\begin{equation} \label{pcrtbpy} \ddot{y} +2\dot{x}=y-(1-\mu)\frac{y}{r_{1}^{3}}-\mu \frac{y}{r_{2}^{3}} \end{equation}
where $r_{1} = \sqrt{(x+\mu)^{2} + y^{2}}$ is the distance from the spacecraft to the primary body and $r_{2} = \sqrt{(x-1+\mu)^{2} + y^{2}} $ is the distance to the secondary. Figure \ref{fig:crtbp} is a diagram of the model, except for in our analysis we restrict ourselves to the case of $z = 0$. 

There are two important properties of equations \eqref{pcrtbpx} and \eqref{pcrtbpy} to note. First of all, the Jacobi integral
\begin{equation} \label{jacobi} C = x^{2} + y^{2}  + 2\left( \frac{1-\mu}{r_{1}}+ \frac{\mu}{r_{2}} \right) - \left( \dot x^{2} + \dot y^{2}  \right)  \end{equation}
is a constant of motion. Furthermore, the equations of motion are in fact Hamiltonian, with $H = -\frac{1}{2}C$. Hence, trajectories in the PCRTBP are restricted to 3-dimensional submanifolds of the state space  satisfying $H(x,y,\dot x, \dot y)=$ constant. 

The second property to note is that the equations of motion have a time-reversal symmetry. Namely, if $(x(t), y(t), t)$ is a solution of Equations \eqref{pcrtbpx} and \eqref{pcrtbpy} for $t > 0$, then $(x(-t), -y(-t), t)$ is a solution for $t < 0$. 

\subsection{Delaunay and Synodic Delaunay Coordinates} \label{coordinates}

The PCRTBP model described above admits a change of coordinates from $(x,y,\dot x, \dot y)$ to action angle coordinates,, which will be required for the first-order Melnikov analysis carried out in Section \ref{MelnikovSection} (all other computations in this study will be done in the synodic cartesian coordinate frame). We summarize Celletti \cite{celletti} here. Consider an inertial reference frame centered at the primary body $m_{1}$, and let $m_{2} =0 $. Recall that the planar two-body problem in this coordinate frame can be expressed in Delaunay coordinates $(L_{0}, G_{0}, \ell_{0}, g_{0})$, which are closely related to the classical orbital elements. For a two-body orbit, angle $\ell_{0}$ is the mean anomaly, angle $g_{0}$ is the longitude of periapsis, and actions $L_{0}$ and $G_{0}$ are related to the semi-major axis $a$ and eccentricity $e$ as follows:

\begin{equation} L_{0} = \sqrt{a} \quad \quad G_{0} = L_{0}\sqrt{1-e^{2}}\end{equation}

Other texts generally write $L_{0} = \sqrt{\mathcal{G}m_{1} a}$; however with our normalized units, $\mathcal{G}m_{1} = 1$ in the 2-body problem. In these coordinates, it can be shown that the evolution of $(L_{0}, G_{0}, \ell_{0}, g_{0})$ is Hamiltonian with Hamiltonian function 
\begin{equation}  \label{twoBodyHamiltonian}  H(L_{0}, G_{0}, \ell_{0}, g_{0}) = -\frac{1}{2L_{0}^{2}} \end{equation}
and satisfies Hamilton's equations of motion
\begin{equation} \frac{dL_{0}}{dt} = -\frac{\partial H}{\partial \ell_{0}}=0 \quad \quad \frac{dG_{0}}{dt} = -\frac{\partial H}{\partial g_{0}}=0    \end{equation}
\begin{equation} \frac{d\ell_{0}}{dt} = \frac{\partial H}{\partial L_{0}} \quad \quad \frac{dg_{0}}{dt} = \frac{\partial H}{\partial G_{0}}    \end{equation}
As expected, the actions are constant along trajectories, while only the angles (in fact, only $\ell_{0}$) vary. If one now introduces a second large body of $\mathcal{G}m_{2} = \mu$, then the system Hamiltonian Equation \eqref{twoBodyHamiltonian} becomes 
\begin{equation}  \label{DelaunayHamiltonian}  H(L_{0}, G_{0}, \ell_{0}, g_{0}) = -\frac{1}{2L_{0}^{2}} + \mu H_{1} (L_{0}, G_{0}, \ell_{0}, g_{0},t)\end{equation}
where the perturbation $H_{1} (L_{0}, G_{0}, \ell_{0}, g_{0},t)$ is 
\begin{equation}  \label{DelaunayH1}  
H_{1}  (L_{0}, G_{0}, \ell_{0}, g_{0},t) =  \frac{r_{1} \cos (\theta - t)}{\rho_{2}^{2}} 
- \frac{1}{\sqrt {\rho_{2}^{2} + r_{1}^{2} -2 \rho_{2} r_{1} \cos(\theta - t)}}   
 \end{equation}
The quantity $\rho_{2}$ is the constant distance from $m_{2}$ to $m_{1}$; with our normalized units, $\rho_{2} = 1$.  $r_{1}$ as defined earlier is the distance from the spacecraft to $m_{1}$. $\theta = g_{0} + f$ is the longitude of the spacecraft, where $f$ is the spacecraft instantaneous true anomaly. Note that $r_{1}$ and $f$ are functions of $L_{0}$, $G_{0}$, and $\ell_{0}$. 

Now, make a time-varying canonical change of variables $(L, G, \ell, g) = (L_{0}, G_{0}, \ell_{0}, g_{0} - t)$; the new variable $g$ is the instantaneous longitude of periapsis of the spacecraft orbit relative to the x-axis of the the synodic cartesian coordinate frame. Then, the Hamiltonian function from equations \eqref{DelaunayHamiltonian} and  \eqref{DelaunayH1} becomes
\begin{equation}  \label{FinalHamiltonian}  H(L, G, \ell, g) = -\frac{1}{2L^{2}} - G + \mu H_{1} (L, G, \ell, g)\end{equation}
\begin{equation}  \label{FinalH1}  H_{1}  (L, G, \ell, g) =  \frac{r_{1} \cos (g+f)}{\rho_{2}^{2}} 
- \frac{1}{\sqrt {\rho_{2}^{2} + r_{1}^{2} -2 \rho_{2} r_{1} \cos(g+f)}}   
 \end{equation}
which is no longer time-varying. We henceforth refer to these new coordinates as synodic Delaunay coordinates. Note that in these coordinates, for $\mu = 0$, the actions $L$ and $G$ are constant on trajectories, but 
\begin{equation} \frac{d\ell}{dt} = \frac{\partial H}{\partial L} = \frac{1}{L^{3}} = a^{-3/2} \quad \quad \frac{dg}{dt} = \frac{\partial H}{\partial G}  =-1  \end{equation}
Since both angles are varying with time, even for $\mu = 0$ ($m_{2}$ infinitesimal, the two-body problem), not all orbits are periodic in these synodic Delaunay coordinates. Only orbits such that $k_{1} a^{-3/2} + k_{2} (-1) = 0$ for some $k_{1}, k_{2} \in \mathbb{Z}$ will be periodic, with period $2 \pi k_{1}$. Note that $a^{-3/2}$ is the mean motion of the spacecraft, and $1$ is the mean motion of $m_2$. Hence, for $\mu = 0$, an orbit in these coordinates is periodic if and only if the mean motions of the spacecraft and $m_{2}$ are rational multiples of each other. This is equivalent to there being $n,m \in \mathbb{Z}$ such that  in the inertial reference frame, the spacecraft makes $n$ revolutions around $m_{1}$ in the time that $m_{2} $ makes $m$ revolutions around $m_{1}$. In the two-body problem ($\mu = 0$), such orbits are defined as $n:m$ resonant periodic orbits.

\section{Persistence of Resonant Periodic Orbits} \label{MelnikovSection}

As described in section \ref{coordinates}, for $\mu = 0$, in synodic Delaunay coordinates, the only periodic orbits are $n:m$ resonant periodic orbits, $n, m \in \mathbb{Z}$. We are now interested in seeing which of these periodic orbits survive the perturbation when $\mu > 0$. For this, the perturbative method of Melnikov \cite{guckholmes} is useful here. Without going into a full derivation, the essential theory is that given a periodic orbit $\bold{x}_{0}(t)$ in the $\mu = 0$ system, we can express solutions of the $\mu$-dependent equations of motion (with initial condition $\bold{x}_{\mu}(0) = \bold{x}_{0}(0)$) as an expansion in powers of $\mu$
\begin{equation} \label{muexpansion}   \bold{x}_{\mu}(t) = \bold{x}_{0}(t) + \mu \bold{x}_{1}(t)+O(\mu^{2}) \end{equation}
where $\bold{x}_{\mu}(t) = \left( L(t, \mu), G(t,\mu), \ell(t, \mu), g(t,\mu) \right) $ . 

Denote the period of $\bold{x}_{0}(t)$ by $T = 2 \pi m$. The main conclusion of the Melnikov theory is that if an initial condition $\bold{x}_{0}(0) = (L_{i}, G_{i}, \ell_{i}, g_{i})$ can be found such that $\bold{x}_{1}(T) = \bold{x}_{1}(0)$ in the perturbative expansion equation \eqref{muexpansion}, then a true periodic orbit can be found near $\bold{x}_{0}(0)$ for $\mu$ small enough. This means that we can expect to be able to  continue the $\mu = 0 $ periodic orbit $\bold{x}_{0}(t)$ into $\mu > 0$. Furthermore, if one fixes $L_{i}$ and $G_{i}$, and also (without loss of generality) sets $\ell_{i} = 0$, it can be shown that $\bold{x}_{1}(T) = \bold{x}_{1}(0)$ if and only if the Melnikov function
\begin{align}  \label{Melnikov}  
M(g_{i}) &= \int_{0}^{2 \pi m} \left(\frac{\partial H_{0}}{\partial \ell} \frac{\partial H_{1}}{\partial L} - \frac{\partial H_{0}}{\partial L} \frac{\partial H_{1}}{\partial \ell}\right) (L_{i},G_{i}, \Omega(L_{i}) t, g_{i} - t) \, dt  \\
&= \int_{0}^{2 \pi m} - \frac{1}{L_{0}^{3}} \frac{\partial H_{1}}{\partial \ell} (L_{i},G_{i}, \frac{1}{L_{i}^{3}} t, g_{i} - t) \, dt
\end{align}
has simple zeros. If one of those zeros is at $g_{i} = g_{i,z}$, then we know that the periodic orbit with initial condition $\bold{x}_{0}(0) = (L_{i}, G_{i}, 0, g_{i,z})$ persists in a perturbed form for $\mu > 0$, albeit with a possibly slightly different period. Hence, one studies the Melnikov function $M(g_{i})$ given in equation \eqref{Melnikov}. Note that in the integral for $M(g_{i})$, the integration of $\frac{\partial H_{1}}{\partial \ell} $ occurs only along the original, unperturbed periodic orbit. 

One property of $M(g_{i})$ is that it is an odd function, $M(g_{i}) = - M(-g_{i})$. To show this, first note that 
\begin{equation} H_{1} (L, G, \ell, g) = H_{1} (L, G, -\ell, -g) \end{equation}
which then implies that
\begin{align}  \label{oddpartial} \begin{split}
\frac{\partial H_{1}}{\partial \ell}  (L, G, \ell, g) &=  \frac{\partial }{\partial \ell}  \left[ H_{1}(L, G, -\ell, -g)  \right] \\ 
&=  -\frac{\partial H_{1}}{\partial \ell}  (L, G, -\ell, -g)
\end{split} \end{align}
Hence, we find that (using $s = -t $ below)
\begin{align}  \begin{split}
M(g_{i}) &=  \int_{0}^{2 \pi m} - \frac{1}{L_{0}^{3}} \frac{\partial H_{1}}{\partial \ell} (L_{i},G_{i}, \frac{1}{L_{i}^{3}} t, g_{i} - t) \, dt \\
&=  \int_{0}^{-2 \pi m}  \frac{1}{L_{0}^{3}} \frac{\partial H_{1}}{\partial \ell} (L_{i},G_{i}, -\frac{1}{L_{i}^{3}} s, g_{i} + s) \, ds \\
&=  \int_{-2 \pi m}^{0}  -\frac{1}{L_{0}^{3}} \frac{\partial H_{1}}{\partial \ell} (L_{i},G_{i}, -\frac{1}{L_{i}^{3}} s, g_{i} + s) \, ds \\
(*)&=  \int_{0}^{2 \pi m}  -\frac{1}{L_{0}^{3}} \frac{\partial H_{1}}{\partial \ell} (L_{i},G_{i}, -\frac{1}{L_{i}^{3}} s, g_{i} + s) \, ds \\
(**)&=  \int_{0}^{2 \pi m}  \frac{1}{L_{0}^{3}} \frac{\partial H_{1}}{\partial \ell} (L_{i},G_{i}, \frac{1}{L_{i}^{3}} s, -g_{i} - s) \, ds \\
&= -M(-g_{i}) 
\end{split} \end{align}
where line $(*)$ is because $(L_{i},G_{i}, -\frac{1}{L_{i}^{3}} s, g_{i} + s)$ is a $2 \pi m$-periodic orbit, and the line $(**)$ follows from Equation \eqref{oddpartial}. Hence, we have proven that $M(g_{i})$ is odd, and therefore has a zero at $g_{i} = 0$. 

\begin{figure}
\centering
\includegraphics[width=0.5\textwidth]{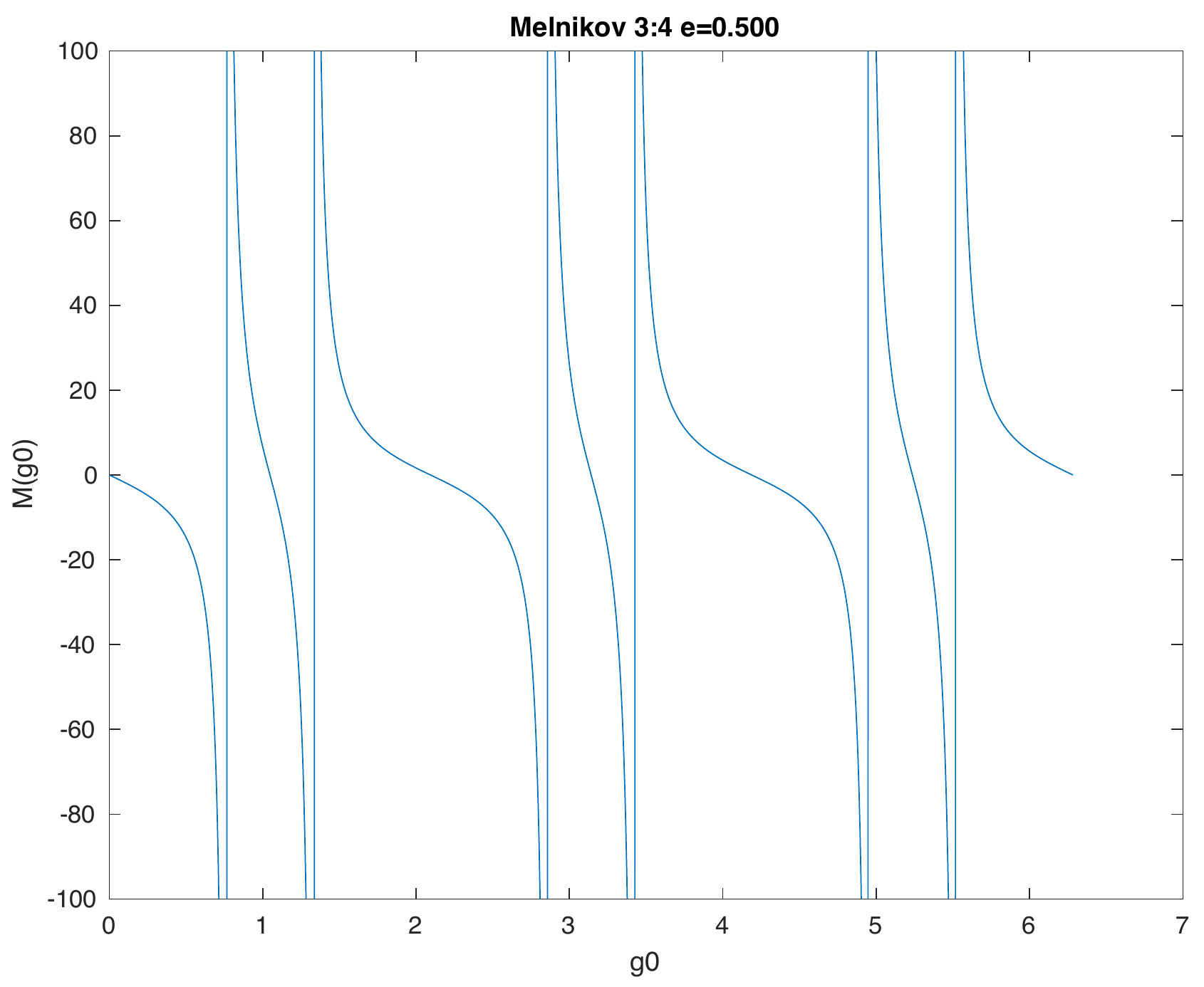}
\caption{\label{fig:melnikovplot}  Plot of $M(g_{i})$ for 3:4 resonance, $e = 0.5$} 
\end{figure}

We plotted $M(g_{i})$ for several different resonances $n:m$. An example of such a plot is shown in Figure \ref{fig:melnikovplot} for $n = 3$, $m=4$, (a 3:4 resonant periodic orbit) with eccentricity $e = 0.5$. 

One thing to note is that $M(g_{i})$ is $2 \pi /n$ periodic when we take $n,m$ coprime. This periodicity is always present, as for an $n:m$ resonant orbit the mean anomaly $\ell = \frac{1}{L_{i}^{3}} t$ is $2 \pi m/n$ periodic. So, evolving the point  $ (L_{i},G_{i}, \ell = 0, g_{i} )$ from $t = 0$ to $t = 2 \pi m/n $ gives the point $ (L_{i},G_{i}, \ell = 0, g_{i} -2 \pi \frac{ m}{n} )$. Both points lie on the same periodic orbit, and so integrating $\frac{\partial H_{1}}{\partial \ell} $ from $t = 0$ to $2\pi m$ along the orbit starting from either point gives the same final result. Integrating starting from the former point corresponds to $M(g_{i})$, while starting from the latter corresponds to $M(g_{i} - 2 \pi \frac{m}{n})$; hence $M(g_{i}) = M(g_{i} - 2 \pi \frac{m}{n})$. Since $n,m$ are coprime, this implies $M(g_{i}) = M(g_{i} + \frac{2 \pi}{n})$. 

However, as is clear from the previous explanation, this periodicity gives us no additional useful zeros of $M(g_{i})$; zeros differing by the quantity $2 \pi /n$ are merely different points on the same orbit, and therefore do not correspond to different persistent resonant orbits. Hence, one can restrict the search for $M(g_{i} ) = 0$ to the interval $g_{i} \in [0, 2 \pi /n)$. Across many different values of $n$ and $m$, apart from $g_{i} = 0$, the only other zero found for all tested cases was $g_{i} = \frac{ \pi}{n}$. We have not analytically proven that $M(\pi/n) = 0$ for arbitrary $m,n$, but the numerical evidence is strong. 

In summary, we have found that the two relevant zeros of $M(g_{i})$ for an $n:m$ resonant periodic orbit are $g_{i} = 0$ and $g_{i} = \frac{\pi}{n}$. Hence, for $\mu > 0$ small enough, it should be possible to find periodic orbits close to the Keplerian orbits with initial conditions $(L_{i}, G_{i}, \ell = 0, g = 0)$ and $(L_{i}, G_{i}, \ell = 0, g = \pi/n)$, where $L_{i} = \sqrt a$ should be chosen so that the corresponding Keplerian orbit period satisfies the $n:m$ resonance condition; $G_{i}$ should satisfy $0 < G_{i} < L_{i}$. Furthermore, as a consequence of the Poincar\'e-Birkhoff fixed point theorem \cite{brown1977}, one of these two orbits will have elliptic stability type and the other should have hyperbolic stability. Intuitively, one expects the orbit corresponding to initial conditions $(L_{i}, G_{i}, \ell = 0, g = 0)$ to be the unstable, hyperbolic orbit, as this corresponds to the initial argument of periapse being aligned with a close flyby of $m_{2}$. It is on resonant orbits of this type that we concentrate now. 

\section{Computation of Resonant Periodic Orbits}

With the persisting Keplerian resonant periodic orbits found, we next compute these surviving orbits and their periods for the full PCRTBP with physically relevant values of $\mu > 0$. Namely, for the Jupiter-Europa system we use $\mu_{E} = 2.5266448850435028 \times 10 ^{-5}$, and for Earth-Moon we used $\mu_{M} = 1.2150584270571545 \times 10^{-2}$. To this end, a  continuation method was used, whereby the periodic orbits computed for smaller values of $\mu$ are used to find an initial guess for the periodic orbit and period corresponding to a larger value of $\mu$. 

We start with a value of $\mu$ for which we wish to compute an $n:m$ resonant orbit. We set $\mu_{0} = 0$, $\mu_1 = \mu/N$, \dots, $\mu_k = k\mu/N$, \dots, $\mu_N = \mu$. We then seek to compute periodic points $\bold{x}_{\mu_k}$ and periods $T_{sc,\mu_k}$ corresponding to the PCRTBP periodic orbit for mass ratio value $\mu_{k}$. $\bold{x}_{\mu_0}$ and $T_{sc,\mu_0} = 2 \pi m$ are known from the Melnikov analysis; to simplify the computations, we convert the initial condition $\bold{x}_{\mu_0} = (L_{i}, G_{i}, \ell = 0, g = 0)$ back to the synodic cartesian coordinate frame $(x_{i},y_{i},\dot x_{i}, \dot y_{i})$ and carry out subsequent computations in that frame. 

To compute the $\bold{x}_{\mu_k}$ and $T_{sc,\mu_k}$, we 
    \begin{enumerate}
    \item Form an initial guess for $(\bold{x}_{\mu_k}, T_{sc,\mu_k}) $ as
    \begin{equation} 
(\bold{x}_{\mu_k}, T_{sc,\mu_k})_{guess} = (\bold{x}_{\mu_{k-1}}, T_{sc,\mu_{k-1}}) + 
 [(\bold{x}_{\mu_{k-1}}, T_{sc,\mu_{k-1}}) - (\bold{x}_{\mu_{k-2}}, T_{sc,\mu_{k-2}}) ]
 \end{equation} except if $k = 1$, $(\bold{x}_{\mu_1}, T_{sc,\mu_1})_{guess} = (\bold{x}_{\mu_0}, T_{sc,\mu_0})$. 
    \item Solve for ($\bold{x}_{\mu_k}$, $T_{sc,\mu_k}$) using initial guess and the MATLAB function fsolve on the equation \begin{equation} \label{fixedpoint}   \Phi_{T_{sc,\mu_k}} (\bold{x}_{\mu_k}) - \bold{x}_{\mu_k} = 0 \end{equation}
    where $\Phi_{T_{sc,\mu_k}}(\bold{x}_{\mu_k})$ denotes the flow of $\bold{x}_{\mu_k}$ by the equations of motion \eqref{pcrtbpx} and \eqref{pcrtbpy} by time $T_{sc,\mu_k}$. 
    \item Increase $k$ by 1, and return to step 1 until $k = N$. 
    \end{enumerate}
Note that $T_{sc,\mu_k}$ must be allowed to vary in order to find periodic orbits for $\mu > 0$. Also, the solution of equation \eqref{fixedpoint} is not unique for a given $\mu_{k}$, as the value of the Jacobi constant (equation \eqref{jacobi}) is not fixed, nor is there a condition added to fix the phasing of the point on a given periodic orbit. Nevertheless, the continuation was successful in continuing 100 different  Keplerian resonant periodic orbits to $\mu = \mu_{E}$, and 32 different orbits to $\mu = \mu_{M}$. We conjecture that for a given resonance $n:m$ (and hence fixed semi-major axis $a$), continuation of orbits with different values of eccentricity $e$ yields final orbits at different values of the Jacobi constant which can be computed from each other through continuation by energy. Additionally, do note that there exist resonant periodic orbits for $\mu > 0$  which are not continuations of $\mu = 0$ orbits \cite{Anderson2016}. 

\begin{figure}
\centering
\includegraphics[width=0.5\textwidth]{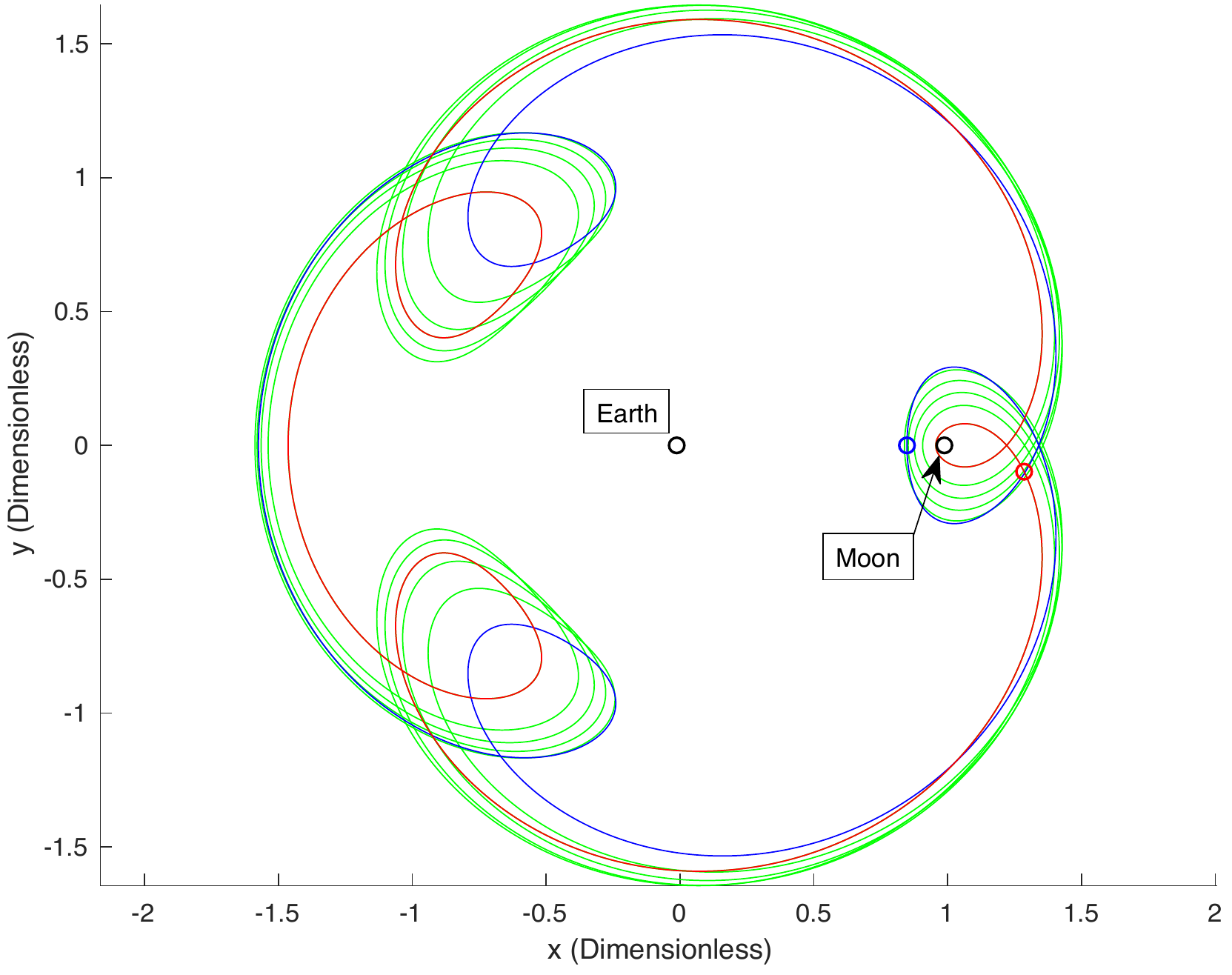}
\caption{\label{fig:continuation}  Continuation of 3:4, $e=0.3, g_{0} =0$ resonant orbit from $\mu = 0$ (blue) to $\mu = \mu_{M}$ (red) with orbits for intermediate $\mu$ values shown in green. } 
\end{figure}

An example of the continuation of a 3:4 resonant orbit with $e = 0.3$ in the Earth-Moon system is shown in Figure \ref{fig:continuation}. The blue curve is the original Keplerian periodic orbit. The red curve is the final computed periodic orbit for $\mu = \mu_{M}$, and the green curves are computed orbits for some intermediate $\mu$ values.  

\section{Parameterization of Invariant Manifolds} \label{parambigsection}

With the resonant periodic orbits and their periods computed for physically relevant values of $\mu$, we next turn our attention to accurate computation of the orbits' stable and unstable invariant manifolds. As mentioned in the introduction, generally current studies using manifolds use linear approximations of invariant manifolds found by computing eigenvectors of the monodromy matrix of the periodic orbit. However, in our case, we compute high order (degree 25 to 50) Taylor polynomials which approximate the manifolds very accurately within some domain of validity. 

Consider a hyperbolic resonant periodic orbit in the PCRTBP containing periodic point $\bold{x}_{\mu}$ and of period $T_{sc,\mu}$. To simplify computations, instead of considering the equations of motion \eqref{pcrtbpx} and \eqref{pcrtbpy}, we instead consider the map $F: \mathbb{R}^{4} \rightarrow \mathbb{R}^{4}$ defined as the time-$T_{sc,\mu}$ mapping by the equations of motion; using the notation established in the previous section, this simply means $F(\bold{x}) = \Phi_{T_{sc, \mu}} (\bold{x})$. 

We know that $\bold{x}_{\mu}$ is a fixed point of the map $F$, and hence the monodromy matrix $DF(\bold{x}_{\mu})$ represents the linearized dynamics around $\bold{x}_{\mu}$. Since we are looking at a hyperbolic periodic orbit, $DF(\bold{x}_{\mu})$ has one stable and one unstable eigenvalue, in addition to  two expected unit eigenvalues. Hence, we know that the stable and unstable manifolds of the fixed point $\bold{x}_{\mu}$ of the full nonlinear map $F$ will also be 1-dimensional. Note that if we consider the full continuous-time flow and the periodic orbit, rather than the map $F$ and its fixed point $\bold{x}_{\mu}$, the stable and unstable manifolds of the periodic orbit are 2-D. Specifically, they are cylinders corresponding to the well-known ``tube dynamics" \cite{RossEtAl}. The stable and unstable manifolds of $\bold{x}_{\mu}$ under $F$ will just be non-closed curves contained on the surface of these cylinders; by integrating points from these curves by the equations of motion, one can compute all the points on the cylindrical manifolds of the periodic orbit. 

Remember that motion in our system is restricted to 3-D submanifolds of the state space corresponding to energy level sets; hence, a given periodic orbit and its stable and unstable manifolds will all be contained within a 3-D submanifold. If we have two periodic orbits at the same energy level, then the 2-D unstable manifold of the first orbit and the 2-D stable manifold of the second orbit will also be contained within a 3-D submanifold. Hence, if the manifolds intersect, they will generically intersect along a curve corresponding to a heteroclinic trajectory. Our final goal is to compute these heteroclinic connections between orbits. 

However, computing 2-D manifolds of periodic orbits and their intersections requires significantly more computational tools and power than for 1-D manifolds of fixed points. Hence, we reduce the dimensionality of our problem through two steps. First of all, we compute 1-D stable and unstable manifolds of the fixed point $\bold{x}_{\mu}$ of the map $F$, rather than 2-D manifolds of orbits. Second, we take a Poincar\'e surface of section (a 2-D submanifold of the 3-D energy submanifold) passing through $\bold{x}_{\mu}$ and compute the 1-D intersection of the 2-D stable and unstable manifolds with the surface of section; this is simply done by propagating points from the 1-D manifolds of the fixed point $\bold{x}_{\mu}$ until their closest intersection with the section. These 1-D intersections of the periodic orbit manifolds with the surface of section simply correspond to stable and unstable manifolds of the fixed point $\bold{x}_{\mu}$ under the Poincar\'e return map. 

\subsection{The Parameterization Method for Invariant Manifolds} \label{paramsection}

The parameterization method is a technique in dynamical systems useful for the computation of several types of invariant geometric structures, including invariant tori as well as stable and unstable manifolds of fixed points, periodic orbits, and tori. It works in both Hamiltonian as well as non-Hamiltonian systems. Haro et al. \cite{haroetal} provide an excellent reference for many applications of this method. The essential idea is that if we have a map $F: M \rightarrow M$ where $M$ is some manifold, and we know that there is an $F$-invariant object diffeomorphic to some model manifold $\mathcal{M}$, then we can solve for an injective immersion $W:\mathcal M \rightarrow M$ and diffeomorphism $f: \mathcal M \rightarrow \mathcal M$ such that the invariance equation
\begin{equation}  \label{invariancequation}   F(W(s)) = W(f(s)) \end{equation}
holds for all $s \in \mathcal M$. We refer to $W$ as the parameterization of the invariant manifold, and $f$ as the internal dynamics on the model manifold $\mathcal M$. Equation \eqref{invariancequation} simply states that $F$ maps the image $W(\mathcal M)$ into itself, so that $W(\mathcal M)$ is the invariant object in the full ambient manifold $M$. 

In our case, we seek to parametrize the 1-dimensional stable and unstable manifolds of the fixed point $\bold{x}_{\mu}$ of $F$. Hence, the ambient manifold $M = \mathbb{R}^{4}$, the model manifold $\mathcal M = \mathbb{R}$, and furthermore we can take $f(s) = \lambda s$, where $\lambda$ is the stable or unstable eigenvalue of $DF(\bold{x}_{\mu})$, depending on which manifold we are trying to compute. Hence, the equation to solve for the parameterization $W(s)$ is 
\begin{equation}  \label{invariancequationfinal}   F(W(s)) - W(\lambda s) =0\end{equation}
where $s \in \mathbb{R}$. We express $W$ as a Taylor series
\begin{equation}  \label{series} W(s) = \bold{x_{\mu}} + \sum_{k \geq 1} W_{k}(s)  \end{equation}
where $W_{k}(s)$ is a monomial of degree $k$ in $s$. The constant term in $W$ is $\bold{x}_{\mu}$, and the linear terms will be the stable or unstable eigenvector of $DF(\bold{x}_{\mu})$ (we take unit length eigenvectors). Hence we need to solve for the higher-order terms $W_{k}(s)$, $k \geq 2$. 

Denote $W_{<k} (s) = \bold{x_{\mu}} + \sum_{j = 1}^{k-1} W_{j}(s)  $.  Assume that we have solved for all $W_{j}(s)$ for $j <k$, so that $F(W_{<k}(s)) - W_{<k}( \lambda s)$ has only $s^{k}$ and higher order terms. Then, the method to solve for $W_{k}(s)$ is:
    \begin{enumerate}
    	\item Find $E_{k}(s)=  [F(W_{<k}(s)) - W_{<k}( \lambda s)]_{k}$, where $[\cdot]_{k}$ denotes the $s^{k}$ term of the RHS. 
	\item Solve for the $s^k$ term $W_{k} (s)$ which when added to $W_{<k} (s)$ cancels $E_k (s)$ in equation \eqref{invariancequationfinal}.
	\begin{align}  \begin{split} \label{correctionwk} -E_{k}(s) &= DF(\bold{x}_{\mu}) W_{k}(s) - W_{k}( \lambda s)   \\
	   &= \left[ DF(\bold{x}_{\mu})  - \lambda^{k} I \right] W_{k}(  s)  \end{split} \end{align}
	\item Set $W_{<k+1} (s) = W_{<k} (s) + W_k(s)$ and return to step 1 until satisfied with the degree of $W$
    \end{enumerate}
We start with $k=2$ and proceed. We elaborate on the computation of the degree $k$ monomial $E_{k}(s)$ from step 1 in section \ref{jettransport}. Equation \eqref{correctionwk} can be derived from the requirement 
\begin{equation}   \left[  F(W_{<k}(s)+W_{k}(s)) - \left(W_{<k}(\lambda s)+W_{k}(\lambda s)\right) \right]_{k}= 0 \end{equation} where as before $[\cdot]_{k}$ denotes taking the $s^{k}$ term of the quantity inside brackets.
Expanding the LHS in Taylor series and discarding terms of polynomial degree greater than $k$ gives
\begin{align} 
[F(W_{<k}(s)) +&DF(W_{<k}(s)) W_{k}(s) - \left(W_{<k}(\lambda s)+W_{k}(\lambda s)\right)]_{k} \\
&=E_{k}(s) + [DF(W_{<k}(s)) W_{k}(s) - W_{k}(\lambda s) ]_{k} \\
&=E_{k}(s) + DF(\bold{x}_{\mu}) W_{k}(s) - W_{k}(\lambda s) = 0 
 \end{align}
where the last line follows from the preceding one because one can divide $s^{k}$ out from $E_{k}(s)$, $W_{k}(s)$, and $W_{k}(\lambda s)$, and then take $s \rightarrow 0$.

\subsection{Computing $E_{k}(s)$: Automatic Differentiation and Jet Transport} \label{jettransport}

In step 1 of the parameterization method algorithm, we computed the quantity 
\begin{equation} E_{k}(s)=  [F(W_{<k}(s)) - W_{<k}( \lambda s)]_{k} \end{equation}
$W_{<k}(s)$ is a degree $k-1$ polynomial and $\lambda$ is a constant, so the degree $k$ term of $W_{<k}(\lambda s)$ is just 0. However, $F$ is the nonlinear time-$T_{sc,\mu}$ mapping of phase space points by the equations of motion \eqref{pcrtbpx} and \eqref{pcrtbpy}; hence, computing $F(W_{<k}(s)) $ as a polynomial is not a trivial matter. For this, the tools of automatic differentiation \cite{haroetal} and  jet transport \cite{perezpalau2015} are useful. 

Automatic differentiation is a technique which allows for rapid recursive evaluation of operations on polynomials. For instance, given two polynomials $f(x)$ and $g(x)$, suppose we wish to compute $d(x) = f(x) / g(x)$ as a polynomial. We know that $d(x) = f(x) / g(x) \iff f(x) = d(x) g(x)$; hence, using subscript $j$ to denote the degree $j$ coefficient,
\begin{align}  \begin{split}
f_{k}(x) &= \sum_{j=0}^{k} d_{j}(x) g_{k-j}(x) \\
&=  \sum_{j=0}^{k-1} d_{j}(x) g_{k-j}(x)+ d_{k}(x) g_{0}(x)  \\
\end{split} \end{align}
\begin{equation}  \label{autodiff} \therefore d_{k}(x) = \frac{1}{g_{0}} \left( f_{k}(x) - \sum_{j=0}^{k-1} d_{j}(x) g_{k-j}(x) \right)\end{equation}  
We know that $d_{0} = f_{0}/ g_{0}$, and using equation \eqref{autodiff} with the known coefficients of $f$ and $g$, can recursively find $d_{k}(x)$, $k \geq 1$. Similar recursive formulas exist for $f(x)^{\alpha}$ as well as many other functions \cite{haroetal}. The key property of all automatic differentiation formulas is that the $s^{k}$ coefficient of the output depends only on the $s^{k}$ and lower order coefficients of the operands. Hence, truncation of Taylor series for the purpose of implementation on a computer does not affect the accuracy of the computed coefficients.

The utility of automatic differentiation is that it allows us to substitute polynomials such as $W_{<k}(s)$ for $(x,y,\dot x, \dot y)$ in the equations of motion \eqref{pcrtbpx} and \eqref{pcrtbpy} to get polynomials in $s$ for $(\dot x, \dot y,\ddot x, \ddot y)$. In particular, let $V(s,t) = \sum_{i=0}^{\infty}V_{i}(t)s^{i}:\mathbb{R}^{2} \rightarrow \mathbb{R}^{4}$ be a Taylor series-valued function of time, with time-varying coefficients $V_{i}(t)$. Denote the $x$, $y$, $\dot x$, and $\dot y$ components of $V(s,t)$ as $V_{x}(s,t)$, $V_{y}(s,t)$, ${V}_{\dot{x}}(s,t)$, and $V_{\dot y}(s,t)$. Substituting $V$ in the equations of motion, we get the system of differential equations 
\begin{equation}  \frac{d}{dt}{V_{x}(s,t)} = V_{\dot x}(s,t) \end{equation}
\begin{equation}  \frac{d}{dt}{V_{y}(s,t)} = V_{\dot y}(s,t) \end{equation}
\begin{align}  \begin{split} \frac{d}{dt}{V_{\dot x}(s,t)} = 2V_{\dot y}&(s,t) + V_{x}(s,t) 
-(1-\mu)\frac{V_{x}(s,t)+\mu}{r_{1}(s,t)^{3}}  -\mu \frac{V_{x}(s,t)-1+\mu}{r_{2}(s,t)^{3}}  \end{split} \end{align}
\begin{align}  \begin{split} \frac{d}{dt} V_{\dot y}(s,t)= -2V_{\dot x}(s,t) + V_{y}(s,t) 
-(1-\mu)\frac{V_{y}(s,t)}{r_{1}(s,t)^{3}} -\mu \frac{V_{y}(s,t)}{r_{2}(s,t)^{3}}  \end{split} \end{align}
where $r_{1}(s,t) = \sqrt{(V_x+\mu)^{2} + V_y^{2}}$ and $r_{2}(s,t) = \sqrt{(V_x-1+\mu)^{2} + V_y^{2}} $. For a given $t$, the RHS of each equation can be simplified to a polynomial using automatic differentiation. Hence, this can be interpreted as a differential equation for the polynomial coefficients of each component of $V(s,t)$; for each equation, one simply sets the time derivative of the $s^{m}$ coefficient on the LHS to the $s^{m}$ coefficient on the RHS.  Solving this equation with initial condition $V(s,0) = W_{<k}(s)$, we have that $V(s,T_{sc, \mu}) = F(W_{<k}(s))$, which is the polynomial we need. 

Hence, by treating the coefficients of $W_{<k}(s)$ as real parameters to be integrated from $0$ to $T_{sc, \mu}$, we can numerically integrate $W_{<k}(s)$ coefficient by coefficient to find $F(W_{<k}(s))$. This method of integrating a polynomial curve is known as jet transport; for more details, see Perez-Palau \cite{perezpalau2015}. The essential idea is to overload algebraic operations and numerical integration routines with the ability to accept arrays of polynomial coefficients rather than only floating point numbers. We can use truncated Taylor series in this algorithm since the automatic differentiation formulas used for the evaluation of time derivatives are valid for truncated series. Note that if we have an $n$-dimensional state ($n = 4$ in our case) and a degree-$d$ truncated series, then the integration required is $n(d + 1)$ dimensional. 

\subsection{Notes About Computation of Manifolds} 

\begin{figure}
\centering
\includegraphics[width=0.5\textwidth]{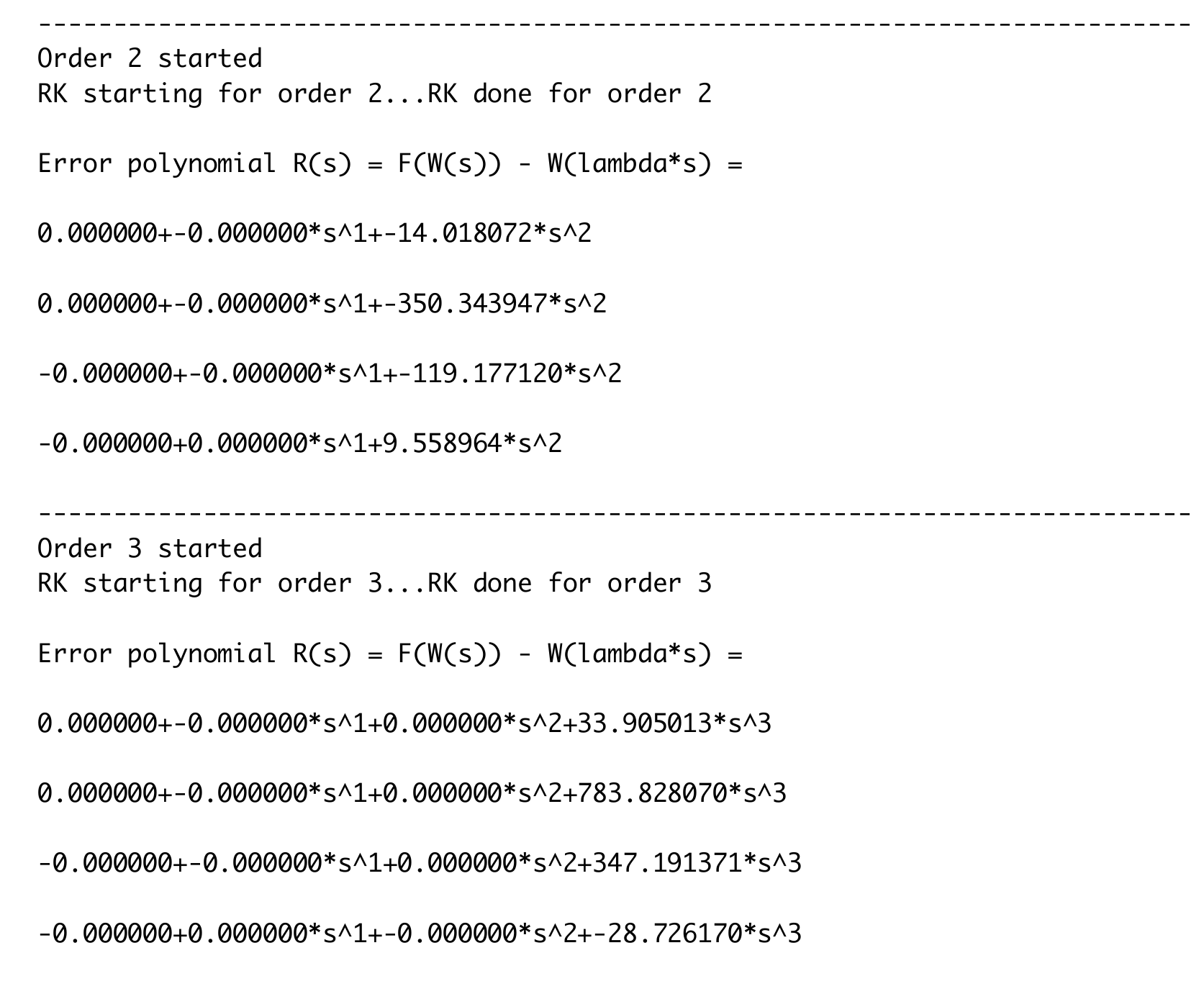}
\caption{  \label{fig:programoutput}Program Output} 
\end{figure}

The parameterization method, automatic differentiation, and jet transport described in the preceding sections were implemented in programs written in C using the GSL library \cite{gsl} for the computation of stable and unstable manifolds. Figure \ref{fig:programoutput} gives an example of part of the program output; in the order $k$ step of the program, first $E_{k}(s) = F(W_{<k} (s)) - W_{<k}(\lambda s) $ is computed using the GSL rk8pd integrator for jet transport (denoted RK in Figure \ref{fig:programoutput}). Printing $E_{k}(s)$ to the terminal, we see that the coefficients of order less than $k$ are zero as expected in each step. The final $d$ degree polynomial $W_{\leq d}(s)$ satisfies $F(W_{\leq d}(s)) - W_{\leq d}(\lambda s) =0$ up to polynomial terms of order $d$. 

To optimize computational time and storage requirements, at the order $k$ step, we only store polynomial coefficients up to degree $k$ in the automatic differentiation and jet transport steps. This allows the jet transport to run much more quickly than it did when we stored additional unnecessary terms. Also, note that if $W(s)$ solves Equation \eqref{invariancequation}, then so does $W(\alpha s)$ where $\alpha$ is an arbitrary constant. Hence, if the jet transport integration is struggling to converge due to fast-growing coefficients of $W(s)$, it helps to scale $W(s)$ to $W(\alpha s)$ by multiplying the eigenvector $W_{1}(s)$ by $\alpha <1$ and then restarting the parameterization method algorithm from Section \ref{paramsection}.

Finally, one last remark is that if one takes the original periodic point $\bold{x}_{\mu}$ to be on the hyperplane $y = 0$, then using the time-reversal symmetry mentioned in Section \ref{modelsection}, we can see that the unstable manifold $W^{u}(s)$ can be found from the stable manifold $W^{s}(s)$ simply by setting $W^{u}(s) = W^{s}(s)$ and then multiplying the $y$ and $\dot x$ components of $W^{u}(s)$ by $-1$. This enables us to save half of the computation time that computing both $W^{s}$ and $W^{u}$ would have taken. Henceforth, we always take $\bold{x}_{\mu}$ on $y = 0$, and  always use this symmetry to compute the unstable manifolds. 

\subsection{Fundamental Domains of Parameterizations} 

Though the $d$ degree polynomial parameterizations $W_{\leq d}(s)$ of the stable and unstable manifolds of $\bold{x}_{\mu}$ are expected to be much more accurate than their linear approximations, they are still inexact and subject to some error. In addition, even if the polynomials could be fully and exactly computed, they still will only be valid within some radius of convergence. Hence, one must determine for which values of $s \in \mathbb{R}$ the polynomial $W_{\leq d}(s)$ is an accurate representation of the invariant manifold. 

To do this, we fix an error tolerance, such as say $E_{tol} = 10^{-5}$ or $10^{-6}$. We then seek to find what is referred to as the fundamental domain of $W_{\leq d}(s)$. The fundamental domain is defined as the maximum magnitude of $s$ such that the error in invariance Equation \eqref{invariancequation} is less than $E_{tol}$. To be precise, we want a  $D \in \mathbb{R}$ such that for all $s$ such that $|s| \leq D $, 
\begin{equation} \|F(W_{\leq d}(s)) - W_{\leq d}( \lambda s)\| < E_{tol} \end{equation}
	
By computing the fundamental domains for over 60 resonant orbit stable manifolds,  we observed orders of magnitude improvement in fundamental domains for $d = 25$ compared to $d = 1$. For linear parameterizations ($d = 1$), the domains of all test cases were on the order of $10^{-4}$ at best, generally $10^{-5}$. However, for the degree-25 polynomial parameterizations $W_{d \leq 25}(s)$, most domains were on the order of 0.1 or even 1. 

Note that if one scales $W_{d \leq 25}(s)$ to $W_{\leq 25}(\alpha s)$ with $\alpha <1$, then the fundamental domain increases by a factor of $\alpha^{-1}$. Hence, whenever we compare domains between parameterizations, we always multiply the domain by any scale factor $\alpha$ used, so that valid comparisons can be made. 

\subsection{Globalization and Visualization} \label{globoviz}

With the fundamental domains computed, we now seek to use the manifold parameterization $W(s)$ to find heteroclinic connections between different resonant periodic orbits. Before we can accomplish this, it is useful to plot the intersection of the periodic orbits' invariant manifolds with a Poincar\'e section. Additionally, we need to compute the manifold $W(s)$ for $s$ values outside the fundamental domain, referred to as globalization. We do these two tasks simultaneously. 

In our case, the Poincar\'e section we use is a $y = 0$, $x<0$ section; recall that the Jacobi constant $C$ is a constant of motion, so fixing the values of $C$ and $y$ restricts us to a 2-D surface of section. We know that the curve $W(s)$ computed earlier is an invariant manifold for the $F$-fixed point $\bold{x}_{\mu}$.  The  curve $W(s)$ lies on the 2-D invariant manifold of the resonant periodic orbit passing through $\bold{x}_{\mu}$. Hence, if we seek to find the intersection of this 2-D invariant manifold with the 2-D surface of section in the 3-D energy level submanifold, this will be a 1-dimensional curve which can be found by propagating points from the curve $W(s)$ to the section. As we are taking $\bold{x}_{\mu}$ to be on $y = 0$, only a short forwards or backwards integration should be required at each point $W(s)$. 

Denote the point found by propagating $W(s)$ to the surface of section as $W_{p}(s)$. Henceforth, denote the forwards and backwards Poincar\'e maps by $P_{+}$ and $P_{-}$, respectively.  Since $F(W(s)) = W(\lambda s)$ (at least within $E_{tol}$), we have that $P_{+} (W_{p}(s)) = W_{p}(\lambda s)$, and that $W_{p}(s)$ is a curve representing the invariant manifold for the fixed point $\bold{x}_{\mu}$ under $P_{+}$. In practice, we take a discrete grid of $s$-values $\{s_{i}\}$ from $-D$ to $D$ (the fundamental domain value), and compute and store $W_{p}(s_{i})$ for each $s_{i}$. For each $W_{p}(s_{i})$, we plot the values $(x, \dot x)$ since given $C$ and $y = 0 $, this is sufficient to determine $\dot y$. 

Note that we no longer have a polynomial representing the manifold on the Poincar\'e section. Instead, we have an accurate grid of points of the manifold $W_{p}(s)$. Computing the polynomial representation of the manifold $W_{p}(s)$ on the Poincar\'e section requires expansion of each coefficient of $W(s)$ as a Taylor series in time under jet transport, followed by the computation of the Poincar\'e return time as a polynomial in $s$ and the composition of the two polynomials, as is described by Perez-Palau \cite{perezthesis}. Rather than carrying out this complicated procedure, we chose to simply propagate a fine grid of points to the section. 

Next, we compute the manifold $W_{p}(s)$ for $s$ values outside the fundamental domain. For this, we now follow the usual process of globalization of invariant manifolds, which is to propagate the fundamental domain \cite{haroetal}. Namely, we take the points $W_{p}(s_{i})$, and propagate them to define $W_{p}(s)$ at larger $s$-values using the equations 
	\begin{equation} W_{p}(\lambda s) = P_{+}(W_{p}(s)) \text{ if } \lambda > 1 \end{equation}
	\begin{equation} W_{p}(s/\lambda) = P_{-}(W_{p}(s)) \text{ if } \lambda < 1 \end{equation}

\begin{figure}
\centering
\includegraphics[width=0.6\textwidth]{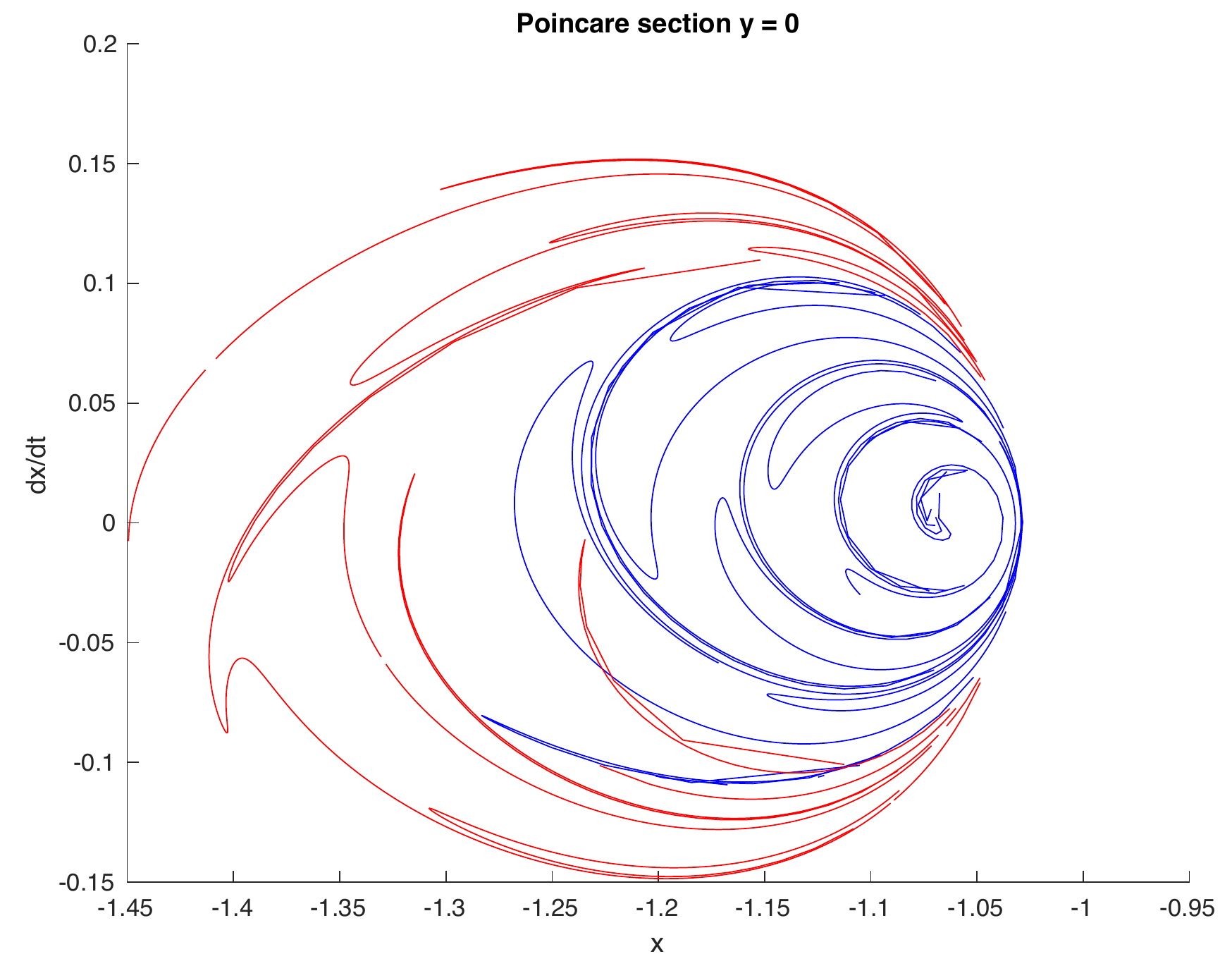}
\caption{ \label{fig:Poincaremap} 3:4 $W^{u}$ (red) and 5:6 $W^{s}$ (blue) Poincar\'e Section for Jacobi Constant $C = 3.0024$}
\end{figure}

We then store the points of $W_{p}(s)$ found and their corresponding $s$-values in a data file. In practice, it is helpful to only count intersections such that $\dot y$ has the same sign as $\dot y$ at $\bold{ x}_{\mu}$. An example Poincar\'e section after globalization, with stable and unstable manifolds of 5:6 and 3:4 resonant orbits, respectively, is given in Figure \ref{fig:Poincaremap}. 

\section{Computation of Heteroclinic Connections } \label{heteroclinic}

With the stable and unstable manifolds of the PCRTBP resonant periodic orbits accurately parametrized, globalized, and plotted on the Poincar\'e section, we now demonstrate how to use the results of the previous computations to find heteroclinic connections between orbits. From now on, denote $W^{u}_{1}(s_{u})$ and $W^{s}_{2}(s_{s})$ as the intersections with the Poincar\'e section of the stable and unstable manifolds of periodic orbits 1 and 2, respectively. Heteroclinic connections from orbit 1 to orbit 2 correspond to intersections of the curves $W^{u}_{1}$ and $W^{s}_{2}$. 

We have the values of $W^{u}_{1}(s_u)$ and $W^{s}_{2}(s_{s})$ on the Poincar\'e section on a discrete grid of $s_{u}$ and $s_{s}$ values, say $\{s_{u,i}\}$ and $\{s_{s,j}\}$. $W^{u}_{1}(s_u)$ and $W^{s}_{2}(s_{s})$ are hence stored as sequences of consecutive points $\{W_{1}^{u}(s_{u,i})\}$ and $\{W_{2}^{s}(s_{s,j})\}$ whose $(x,\dot x)$ values are plotted in the Poincar\'e section. We seek to find $s_{u}$ and $s_{s}$ such that $W^{u}_{1}(s_u) = W^{s}_{2}(s_{s})$. 
To accomplish this numerically, the first part of the algorithm is to:
\begin{figure}
\includegraphics[width=0.485\textwidth]{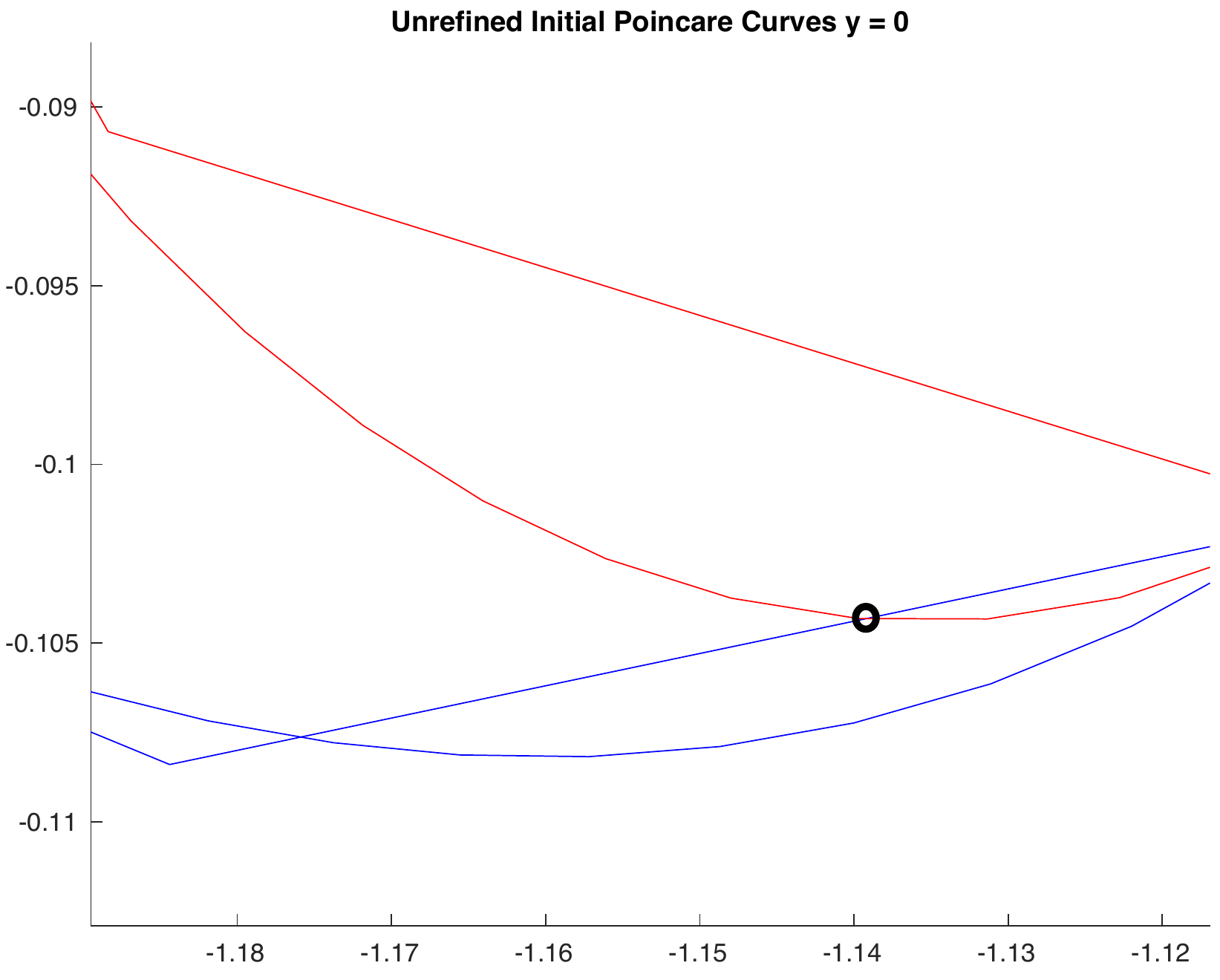}
\includegraphics[width=0.485\textwidth]{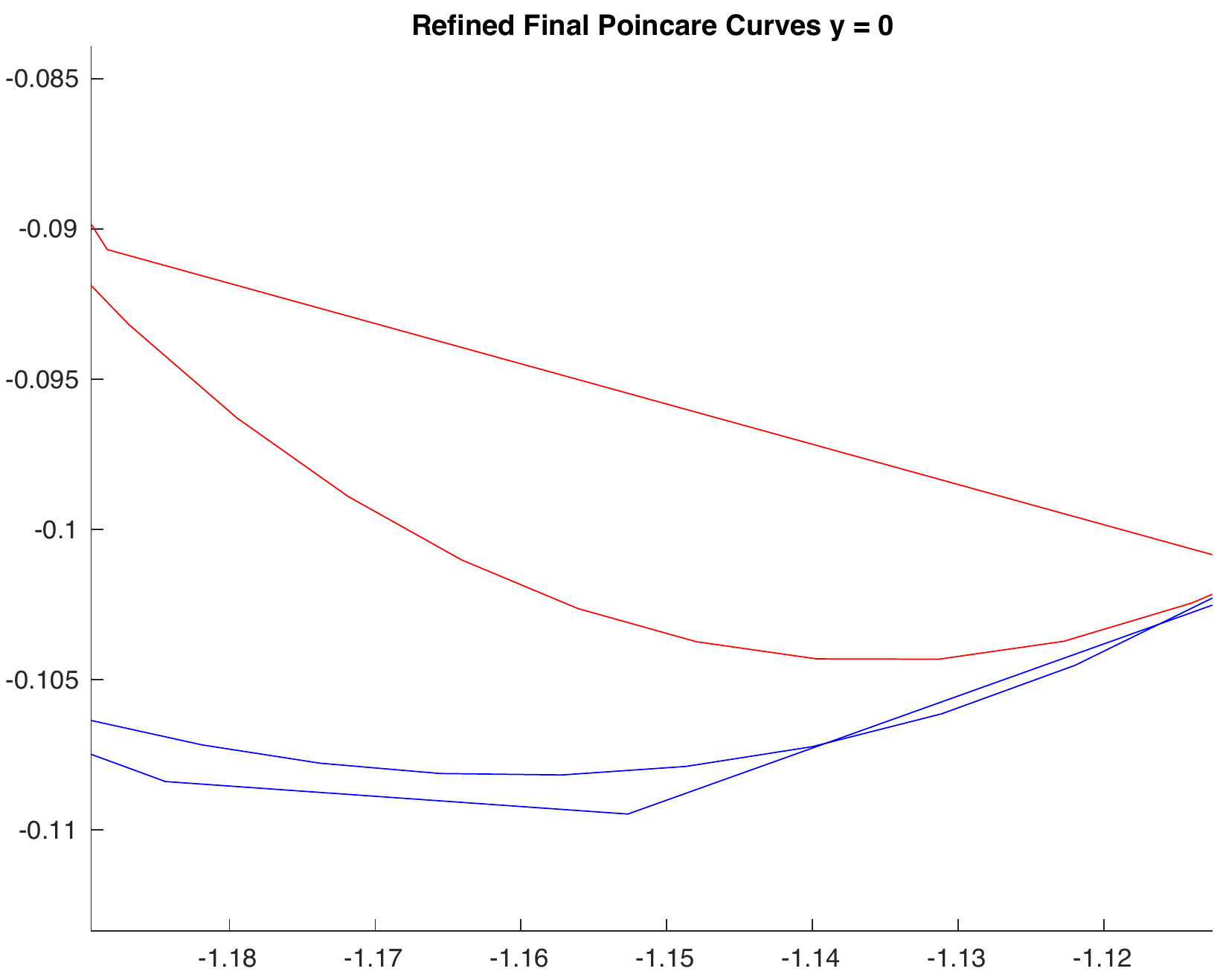}
\caption{ \label{fig:brokenintersect} False intersection (circled) removed upon refinement}
\end{figure}
    \begin{enumerate}
    	\item Connect all consecutive $(x, \dot x)$ points $W_{1}^{u}(s_{u,i})$ and $W_{1}^{u}(s_{u,i+1})$ by line segments (similarly for all $W_{2}^{s}(s_{s,j})$ and $W_{2}^{s}(s_{s,j+1})$)
	\item Remove all line segments corresponding to discontinuities. 
	\item For each segment between points of $W_{1}^{u}$ check for intersections with all segments of $W_{2}^{s}$ 
   \end{enumerate}
   Step 2 is somewhat heuristic; to detect discontinuities, we checked if the quantity $ W_{1}^{u}(s_{u,i+1}) - W_{1}^{u}(s_{u,i}) $ had large values, or if it was much larger in magnitude than $W_{1}^{u}(s_{u,i}) - W_{1}^{u}(s_{u,i-1})$ (similar for $W_{2}^{s}$). Also note that step 3 is easily parallelizable, and indeed benefits significantly from doing so. 

With the first part of the algorithm serving to find intersecting segments of points from $W_{1}^{u}$ and $W_{2}^{s}$, as well as the $s_{u}$ and $s_{s}$ values corresponding to the endpoints, the next part of the algorithm refines the estimate for $s_{u}$ and $s_{s}$ satisfying $W^{u}_{1}(s_u) = W^{s}_{2}(s_{s})$. In particular, if an intersection is detected between the $\{W_{1}^{u}(a_{1}), W_{1}^{u}(b_{1})\}$ segment and $\{W_{2}^{s}(a_{2}), W_{2}^{s}(b_{2})\}$ segment:

	\begin{enumerate}
	\item Find $W_{1}^{u}(\frac{a_{1}+ b_{1}}{2}) = P_{+}^{k}(W_{1}^{u}(\lambda_{u}^{-k}\frac{a_{1}+ b_{1}}{2}))$ where $k$ is such that $\lambda_{u}^{-k}\frac{a_{1}+ b_{1}}{2}$ is in the fundamental domain of the polynomial from which $W_{1}^{u}$ was computed
	\item Find $W_{2}^{s}(\frac{a_{2}+ b_{2}}{2}) = P_{-}^{m}(W_{2}^{s}(\lambda_{s}^{m}\frac{a_{2}+ b_{2}}{2}))$ where $m$ is such that $\lambda_{s}^{m}\frac{a_{2}+ b_{2}}{2}$  is in the fundamental domain of the polynomial from which $W_{2}^{s}$ was computed
	\item Form the segments $\{W_{1}^{u}(a_{1}), W_{1}^{u}(\frac{a_{1}+ b_{1}}{2})\}$, $\{ W_{1}^{u}(\frac{a_{1}+ b_{1}}{2}), W_{1}^{u}(b_{1})\}$  and $\{W_{2}^{s}(a_{2}), W_{2}^{s}(\frac{a_{2}+ b_{2}}{2})\}$, $\{ W_{2}^{s}(\frac{a_{2}+ b_{2}}{2}), W_{2}^{s}(b_{2})\}$
	\item Check for intersections between new segments. If found, return to step 1 with new segment endpoints replacing old ones.
	\item If no intersection found, check for intersections between new segments and other segments on the same continuous curves in $W_{1}^{u}$ and $W_{2}^{s}$. If found, return to step 1. 	
   \item End bisection when $|a_{1}-b_{1}|$ and $|a_{2}-b_{2}|$ are small enough. 
   	\end{enumerate}

In steps 1 and 2, recall that $W_{1}^{u}(\lambda_{u}^{-k}\frac{a_{1}+ b_{1}}{2})$ and $W_{2}^{s}(\lambda_{s}^{m}\frac{a_{2}+ b_{2}}{2})$ are not given directly by the polynomials computed using the parameterization method; however, they are found by integrating points from the polynomials a short distance to the surface of section. Step 5 is necessary because sometimes, when the segments are refined into two segments, intersections that previously existed can break. An example of how this can occur is shown in Figure \ref{fig:brokenintersect}, where the manifolds shown are the same $C = 3.0024$ 3:4 $W^{u}$ and 5:6 $W^{s}$  from Figure \ref{fig:Poincaremap}.

\section{ Example Application to Resonance Transfer in the Jupiter-Europa System}
The methodology described in previous sections is general, and can be applied to systems with a variety of mass ratios $\mu$. In particular, we successfully applied the parameterization method, automatic differentiation, and jet transport to the computation of Taylor series expansions of manifolds in both the Earth-Moon and Jupiter-Europa PCRTBP systems. For the computation of heteroclinic connections, however, we focused our efforts on the Jupiter-Europa system due to the variety of missions currently being planned for that system, such as Europa Clipper \cite{clipper}, Europa Lander \cite{Anderson2019}, and Jupiter Icy Moons Explorer \cite{GRASSET20131}.  

\begin{table}
\centering
\begin{tabular}{rlllll}
\hline
& $5:6$ & $3:4$  \\
\hline
$x$	& -1.231240907544348 & -1.391929713356257  \\
$y$   & 0.000000000000000 & 1.4178538082815e-18			\\
$\dot x$ & 0.000000000000000 & -2.9260154691618e-14		\\
$\dot y$ & 0.371411618064504		& 0.609863420586548		\\
$T_{sc}$ & 38.328135171743014		& 25.338526603095760		\\
$\lambda_{s}$ 			& 0.001256465177783		& 0.011341070996024		\\
$\lambda_{u} $ 			& 795.8835769446018		& 88.175093899915780		\\
\hline
\end{tabular}
\caption{\label{table:X}Initial conditions, periods, and eigenvalues for 3:4 and 5:6 resonant periodic orbits at $C = 3.0024$.}
\end{table}

We used the tools developed in the previous sections for the computation of a 3:4 to 5:6 resonance transfer trajectory in the PCRTBP at the Jacobi constant value 3.0024. The initial conditions, periods $T_{sc}$, and monodromy matrix eigenvalues corresponding to each periodic orbit are given in Table \ref{table:X}. 

Using the parameterization method described in Section \ref{paramsection}, we obtained degree 50 Taylor polynomial expansions representing the stable manifolds of the points in Table \ref{table:X}  under the time $T_{sc}$ map by the equations of motion. By the PCRTBP time-reversal symmetry, we also obtain the unstable manifolds. Next, upon computation of the fundamental domains of these polynomials (using $E_{tol} = 10^{-5}$), we found that the domain for the 5:6 orbit polynomial was approximately 0.9904, while that for the 3:4 orbit was approximately 0.7146. The globalization and Poincar\'e section visualization routine described in Section \ref{globoviz} was then applied to the computed polynomials. Globalization is necessary as the manifold parameterizations, when propagated to the Poincar\'e section, do not intersect within the fundamental domain values of the parameters.

As before, we denote $W_{3:4}^{u}$ and $W_{5:6}^{s}$ as being the unstable and stable manifolds of the  Poincar\'e map fixed points corresponding to the 3:4 and 5:6 orbit points from Table \ref{table:X}. The computed Poincar\'e section with $W_{3:4}^{u}$ and $W_{5:6}^{s}$ was shown earlier in Figure \ref{fig:Poincaremap}. With the Poincar\'e section points computed and stored for both $W_{3:4}^{u}$ and $W_{5:6}^{s}$, we then proceeded to compute heteroclinic connections using the bisection method described in Section \ref{heteroclinic}.

\begin{figure*}[p]
\includegraphics[width=0.5\textwidth]{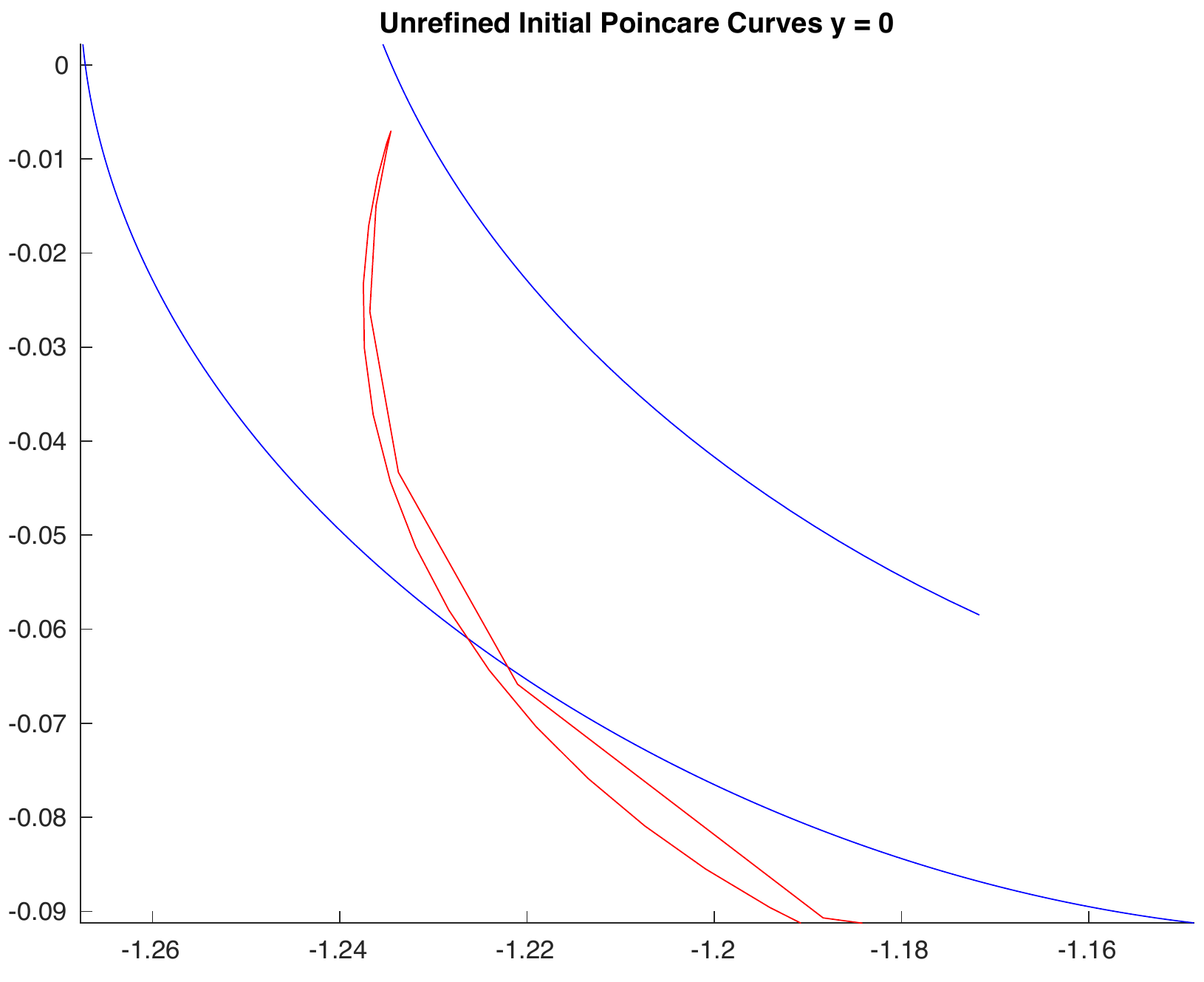}
\includegraphics[width=0.5\textwidth]{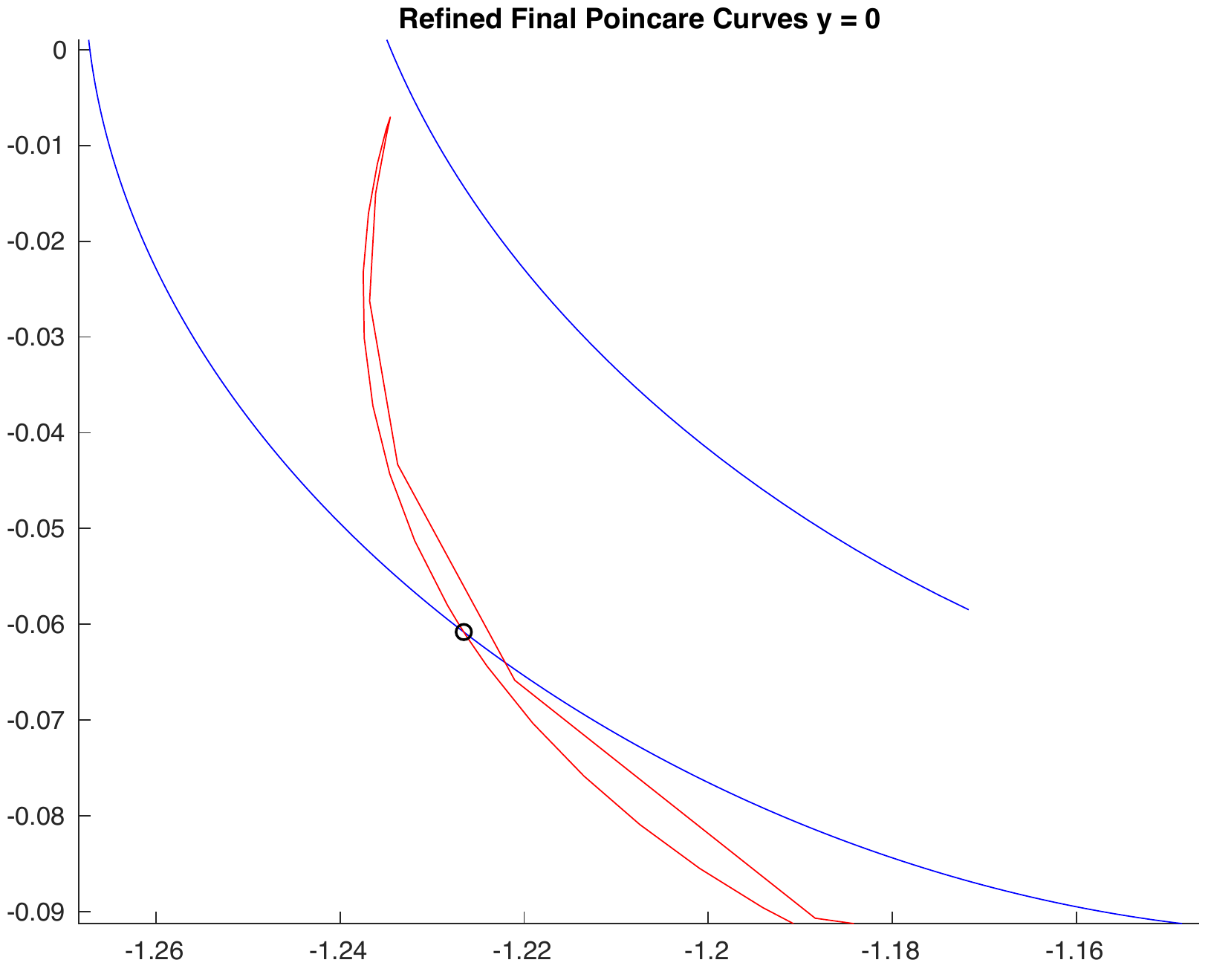}
\includegraphics[width=0.5\textwidth]{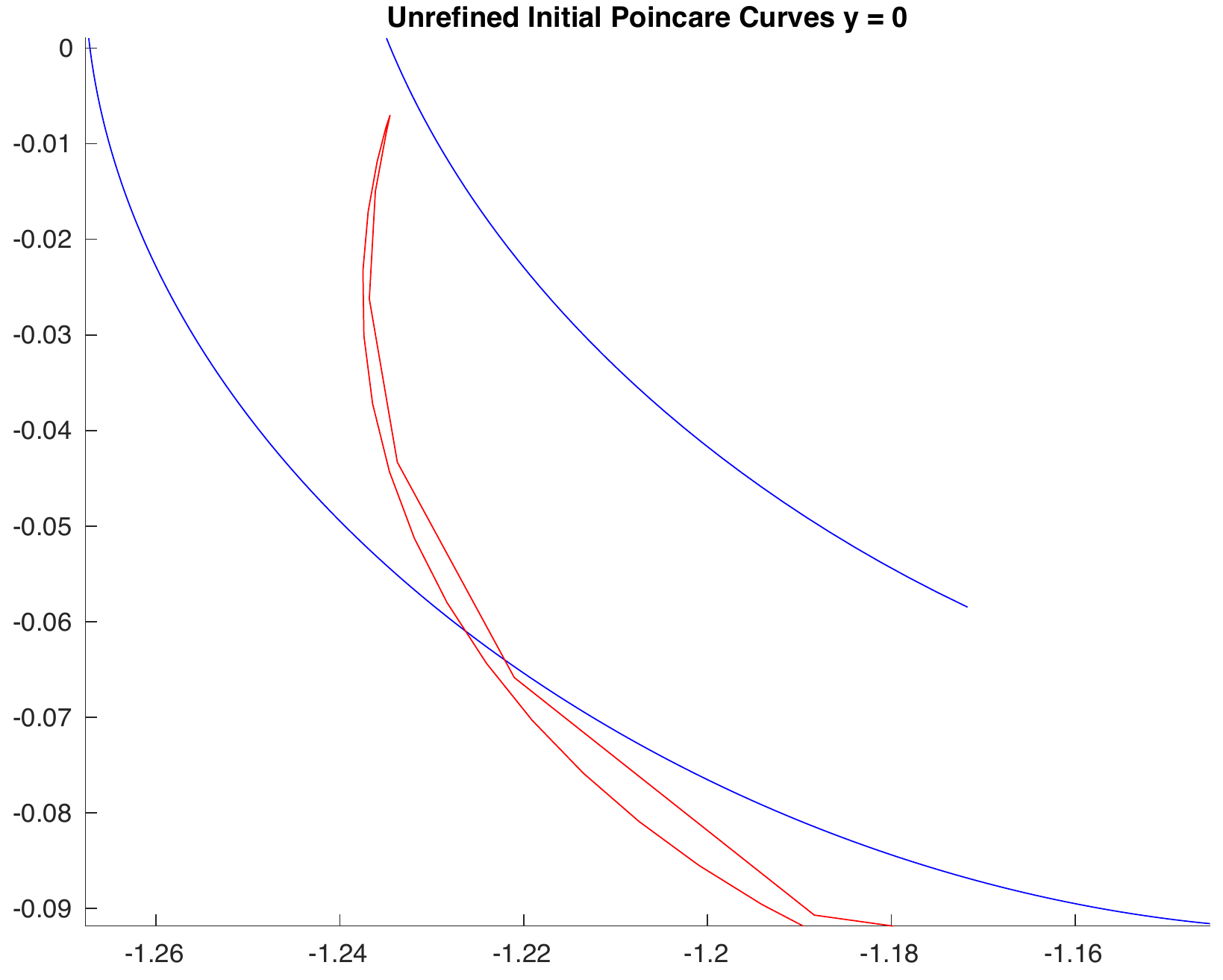}
\includegraphics[width=0.5\textwidth]{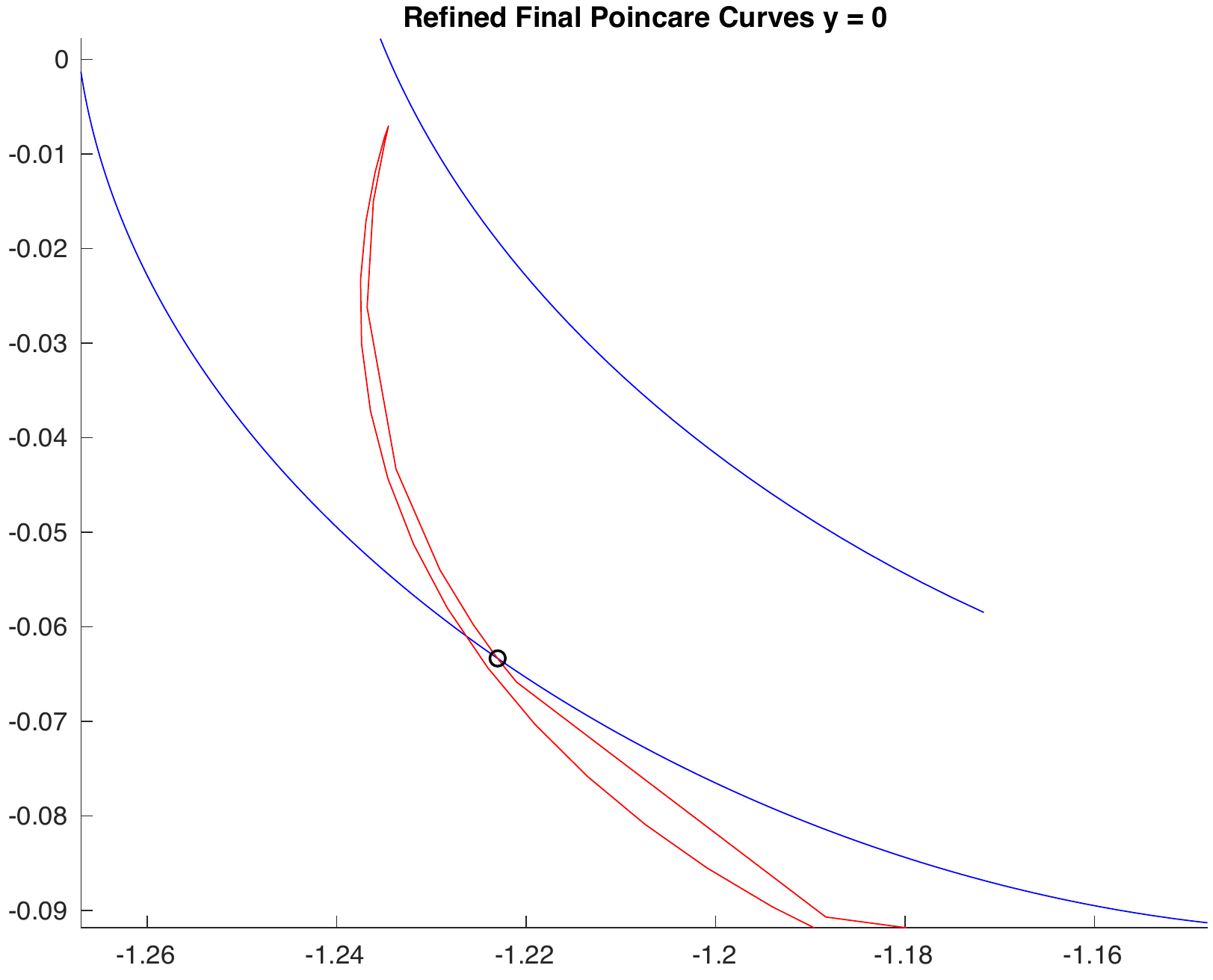}
\includegraphics[width=0.5\textwidth]{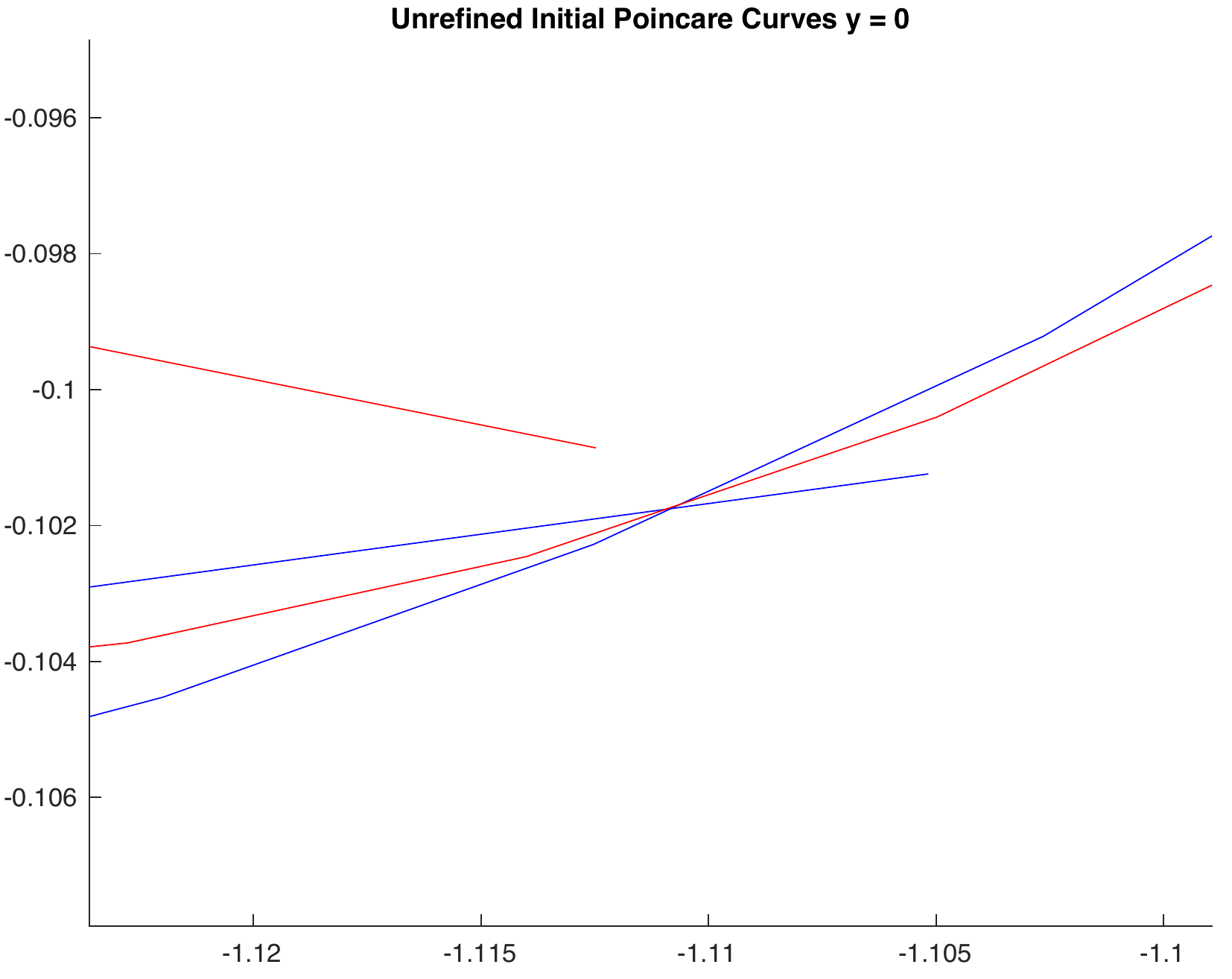}
\includegraphics[width=0.5\textwidth]{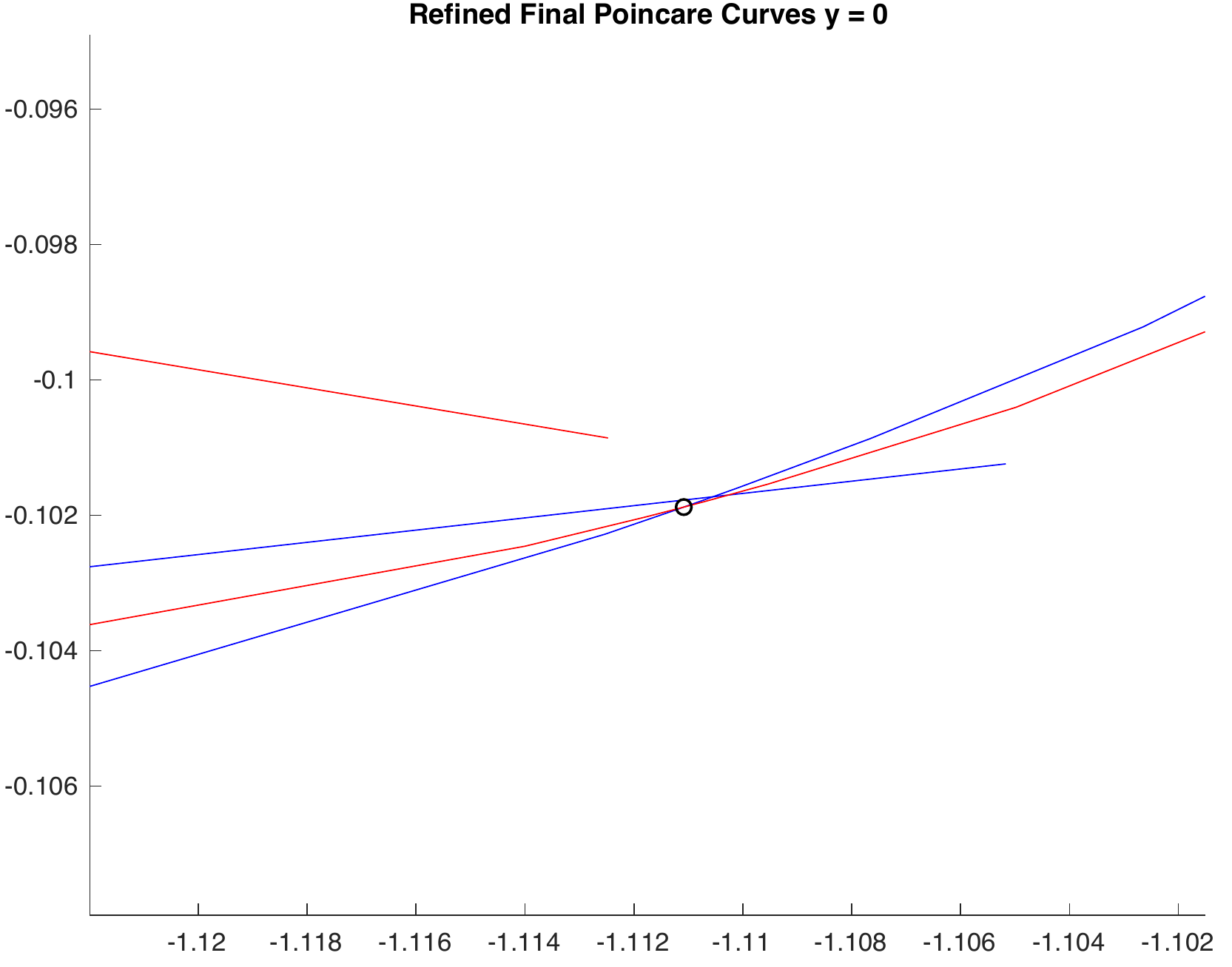}
\caption{ \label{fig:connections} Examples of Approximate Intersections (Left) and Computed Heteroclinic Connections (Right, Circled) for 3:4 to 5:6 Resonance Transfer at Jacobi Constant $C = 3.0024$ ($W_{3:4}^{u}$ red and $W_{5:6}^{s}$ blue)}
\end{figure*}

\begin{figure}
\centering
\includegraphics[width=0.6\textwidth]{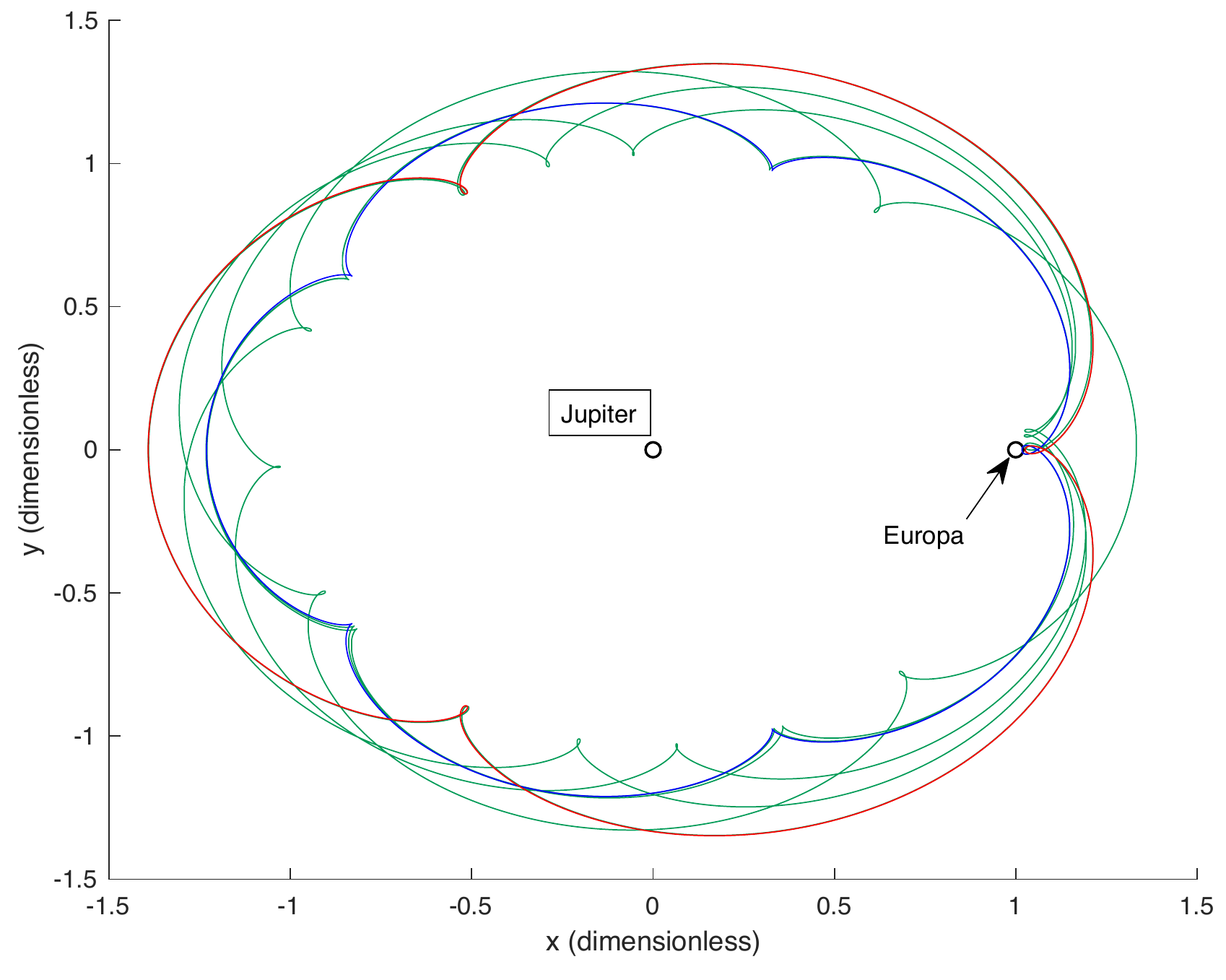}
\caption{ \label{fig:heteroTraj} Trajectory Corresponding to Heteroclinic Connection 3 from Table \ref{table:connections}}
\end{figure}

6 intersections between segments of consecutive stored $W_{3:4}^{u}$ and $W_{5:6}^{s}$ points were initially detected; however, upon refining the segments through the algorithm from Section \ref{heteroclinic}, 3 preliminary intersections were found to be spurious. The coordinates of the 3 computed actual connections are given in Table \ref{table:connections}. Figure \ref{fig:connections} shows how the program refined the Poincar\'e section in the neighborhood of each computed intersection in order to precisely compute the heteroclinic connection point. Finally, Figure \ref{fig:heteroTraj} shows the trajectory corresponding to the third heteroclinic connection point from Table \ref{table:connections}, with the start 3:4 periodic orbit shown in red and the destination 5:6 periodic orbit shown in blue. 

\begin{table}
\centering
\begin{tabular}{rlllll}
\hline
& 1 & 2 & 3  \\
\hline
$x$	& -1.2265598     & -1.2230160     & -1.1110838    \\
$y$   & -4.101840e-14 & -1.989706e-14 & 5.780044e-15 \\
$\dot x$ & -0.060806259   & -0.063340619   & -0.10187786   \\
$\dot y$ & 0.35908692     & 0.35309042     & 0.14762036 \\
$s_{s}$ & -301.609248   & -295.877551 & 14.24735921   \\
$s_{u} $ & -3785.98948   & -3706.35853 & -3874.28227    \\
\hline
\end{tabular}
\caption{\label{table:connections}Computed Heteroclinic Connection Points and corresponding $s_{s}$, $s_{u}$ Values}
\end{table}

Note that our approach of using high order parameterizations of invariant manifolds to compute heteroclinic connections bears some similarity with prior studies; for instance, James and Murray \cite{chebTaylor} parameterized manifolds of periodic orbits using high order Chebyeshev-Taylor series, using the resulting 2D parameterizations to find connecting orbits. However, our study avoids dealing with 2D manifolds by using a Poincar\'e section to reduce the dimensionality of the problem, without sacrificing the accuracy which comes from using high order manifold expansions.

\section{Conclusions}

In this paper, we studied the persistence of resonant periodic orbits in the PCRTBP, and subsequently demonstrated the application of the parameterization method for the computation of high-order expansions of resonant orbit invariant manifolds. We also then demonstrated how to use the resulting polynomials to calculate useful heteroclinic connections. We were able to develop tools to find polynomial approximations of resonant orbit stable and unstable manifolds of degree 25 or even higher; these expansions resulted in a 1000x improvement in the domains of accuracy of the manifold representations as compared to just using linear approximations. 

The tools developed were tested in the Jupiter-Europa system, with the calculations of the manifolds and connections taking only a few minutes on a laptop for a given pair of resonances. The manifold polynomials were used to successfully compute several connections corresponding to 3:4 to 5:6 resonance transition, demonstrating the usefulness of these parameterizations for mission design.

\section*{Acknowledgements}

From June-August 2018 this work was supported by a JPL Strategic University Research Partnership (SURP) grant awarded to Prof. Rafael de la Llave and Dr. Rodney Anderson for FY2018. Part of the research presented in this paper has been carried out at the Jet Propulsion Laboratory, California Institute of Technology, under a contract with the National Aeronautics and Space Administration. From August 2018 onwards, this work was supported by a NASA Space Technology Research Fellowship. The code for the parameterization method was partially adapted from the code developed by Lei Zhang and used in \cite{zhang}. Part of this material is based upon work supported by the National Science Foundation under Grant No. DMS-1440140 while the first and second authors were in residence at the Mathematical Sciences Research Institute in Berkeley, California, during the Fall 2018 semester. The third author is partially supported by NSF grant DMS-1800241.  Paper presented at the 70th International Astronautical Congress, 21-25 October 2019, Washington, D.C., USA; www.iafastro.org. The first author is grateful for support from NASA and the International Space Education Board to attend the IAC 2019. 

\bibliography{bibfile}

\begin{thebibliography}{10}
\expandafter\ifx\csname url\endcsname\relax
  \def\url#1{\texttt{#1}}\fi
\expandafter\ifx\csname urlprefix\endcsname\relax\def\urlprefix{URL }\fi
\expandafter\ifx\csname href\endcsname\relax
  \def\href#1#2{#2} \def\path#1{#1}\fi

\bibitem{Anderson2010}
R.~L. Anderson, M.~W. Lo, \href{https://doi.org/10.2514/1.45060}{Dynamical
  systems analysis of planetary flybys and approach: Planar {Europa} orbiter},
  Journal of Guidance, Control, and Dynamics 33~(6) (2010) 1899--1912.
\newblock \href {http://arxiv.org/abs/https://doi.org/10.2514/1.45060}
  {\path{arXiv:https://doi.org/10.2514/1.45060}}, \href
  {http://dx.doi.org/10.2514/1.45060} {\path{doi:10.2514/1.45060}}.
\newline\urlprefix\url{https://doi.org/10.2514/1.45060}

\bibitem{Anderson2011}
R.~Anderson, M.~Lo, A dynamical systems analysis of resonant flybys: Ballistic
  case, The Journal of the Astronautical Sciences 58.
\newblock \href {http://dx.doi.org/10.1007/BF03321164}
  {\path{doi:10.1007/BF03321164}}.

\bibitem{vaqueroHowell}
M.~Vaquero, K.~C. Howell, \href{https://doi.org/10.2514/1.62230}{Leveraging
  resonant-orbit manifolds to design transfers between libration-point orbits},
  Journal of Guidance, Control, and Dynamics 37~(4) (2014) 1143--1157.
\newblock \href {http://arxiv.org/abs/https://doi.org/10.2514/1.62230}
  {\path{arXiv:https://doi.org/10.2514/1.62230}}, \href
  {http://dx.doi.org/10.2514/1.62230} {\path{doi:10.2514/1.62230}}.
\newline\urlprefix\url{https://doi.org/10.2514/1.62230}

\bibitem{vaqueroCassini}
M.~Vaquero, Y.~Hahn, P.~Stumpf, P.~N. Valerino, S.~V. Wagner, M.~Wong,
  \href{https://arc.aiaa.org/doi/abs/10.2514/6.2014-4348}{Cassini maneuver
  experience for the fourth year of the solstice mission}, in: AIAA/AAS
  Astrodynamics Specialist Conference, 2014.
\newblock \href
  {http://arxiv.org/abs/https://arc.aiaa.org/doi/pdf/10.2514/6.2014-4348}
  {\path{arXiv:https://arc.aiaa.org/doi/pdf/10.2514/6.2014-4348}}, \href
  {http://dx.doi.org/10.2514/6.2014-4348} {\path{doi:10.2514/6.2014-4348}}.
\newline\urlprefix\url{https://arc.aiaa.org/doi/abs/10.2514/6.2014-4348}

\bibitem{Anderson2019}
R.~L. Anderson, S.~Campagnola, D.~Koh, T.~P. McElrath, R.~M. Woollands,
  {Endgame} {Design} for {Europa} {Lander}: {Ganymede} to {Europa} {Approach},
  in: AAS/AIAA Astrodynamics Specialist Conference, no. AAS 19-745, Portland,
  ME, 2019.

\bibitem{Anderson2016}
R.~L. Anderson, S.~Campagnola, G.~Lantoine, Broad search for unstable resonant
  orbits in the planar circular restricted three-body problem, Celestial
  Mechanics and Dynamical Astronomy 124~(2) (2016) 177--199.

\bibitem{guckholmes}
J.~Guckenheimer, P.~Holmes, Nonlinear Oscillations, Dynamical Systems, and
  Bifurcations of Vector Fields, Vol.~42 of Applied Mathematical Sciences,
  Springer, New York, NY, 1983.

\bibitem{CabreFontichLlave}
X.~Cabr{\'e}, E.~Fontich, R.~de~la Llave,
  \href{http://www.sciencedirect.com/science/article/pii/S0022039604005170}{The
  parameterization method for invariant manifolds iii: overview and
  applications}, Journal of Differential Equations 218~(2) (2005) 444 -- 515.
\newblock \href {http://dx.doi.org/https://doi.org/10.1016/j.jde.2004.12.003}
  {\path{doi:https://doi.org/10.1016/j.jde.2004.12.003}}.
\newline\urlprefix\url{http://www.sciencedirect.com/science/article/pii/S0022039604005170}

\bibitem{haroetal}
{\`A}.~Haro, M.~Canadell, J.~Figueras, A.~Luque, J.~Mondelo, The
  Parameterization Method for Invariant Manifolds: From Rigorous Results to
  Effective Computations, Vol. 195 of Applied Mathematical Sciences, Springer
  International Publishing, 2016.

\bibitem{celletti}
A.~Celletti, Stability and Chaos in Celestial Mechanics, Springer-Verlag,
  Berlin; published in association with Praxis Publishing, Chichester, 2010.
\newblock \href {http://dx.doi.org/10.1007/978-3-540-85146-2}
  {\path{doi:10.1007/978-3-540-85146-2}}.

\bibitem{ccar}
\href{https://web.archive.org/web/20181204110444/http://ccar.colorado.edu/geryon/crtbp.html}{The
  circular-restricted three-body problem {(CRTBP)}. from
  https://web.archive.org/web/20181204110444/
  http://ccar.colorado.edu/geryon/crtbp.html} (Retrieved October 6 2020).
\newline\urlprefix\url{https://web.archive.org/web/20181204110444/http://ccar.colorado.edu/geryon/crtbp.html}

\bibitem{brown1977}
M.~Brown, W.~D. Neumann, \href{https://doi.org/10.1307/mmj/1029001816}{Proof of
  the poincar{\'e}-birkhoff fixed point theorem.}, Michigan Math. J. 24~(1)
  (1977) 21--31.
\newblock \href {http://dx.doi.org/10.1307/mmj/1029001816}
  {\path{doi:10.1307/mmj/1029001816}}.
\newline\urlprefix\url{https://doi.org/10.1307/mmj/1029001816}

\bibitem{RossEtAl}
W.~S. Koon, M.~W. Lo, J.~E. Marsden, S.~D. Ross,
  \href{https://doi.org/10.1063/1.166509}{Heteroclinic connections between
  periodic orbits and resonance transitions in celestial mechanics}, Chaos: An
  Interdisciplinary Journal of Nonlinear Science 10~(2) (2000) 427--469.
\newblock \href {http://arxiv.org/abs/https://doi.org/10.1063/1.166509}
  {\path{arXiv:https://doi.org/10.1063/1.166509}}, \href
  {http://dx.doi.org/10.1063/1.166509} {\path{doi:10.1063/1.166509}}.
\newline\urlprefix\url{https://doi.org/10.1063/1.166509}

\bibitem{perezpalau2015}
D.~P{\'e}rez-Palau, J.~J. Masdemont, G.~G{\'o}mez,
  \href{https://doi.org/10.1007/s10569-015-9634-3}{Tools to detect structures
  in dynamical systems using jet transport}, Celestial Mechanics and Dynamical
  Astronomy 123~(3) (2015) 239--262.
\newblock \href {http://dx.doi.org/10.1007/s10569-015-9634-3}
  {\path{doi:10.1007/s10569-015-9634-3}}.
\newline\urlprefix\url{https://doi.org/10.1007/s10569-015-9634-3}

\bibitem{gsl}
\href{https://www.gnu.org/software/gsl/doc/html/}{{GNU} {Scientific} {Library}.
  from https://www.gnu.org/software/gsl/doc/html/} (Retrieved October 2 2019).
\newline\urlprefix\url{https://www.gnu.org/software/gsl/doc/html/}

\bibitem{perezthesis}
D.~P{\'e}rez~Palau, et~al., Dynamical transport mechanisms in celestial
  mechanics and astrodynamics problems, Ph.D. thesis, Universitat de Barcelona
  (2015).

\bibitem{clipper}
T.~{Bayer}, B.~{Cooke}, I.~{Gontijo}, K.~{Kirby}, Europa clipper mission: the
  habitability of an icy moon, in: 2015 IEEE Aerospace Conference, 2015, pp.
  1--12.

\bibitem{GRASSET20131}
O.~Grasset, M.~Dougherty, A.~Coustenis, E.~Bunce, C.~Erd, D.~Titov, M.~Blanc,
  A.~Coates, P.~Drossart, L.~Fletcher, H.~Hussmann, R.~Jaumann, N.~Krupp, J.-P.
  Lebreton, O.~Prieto-Ballesteros, P.~Tortora, F.~Tosi, T.~{Van Hoolst},
  \href{http://www.sciencedirect.com/science/article/pii/S0032063312003777}{Jupiter
  icy moons explorer (juice): An esa mission to orbit ganymede and to
  characterise the jupiter system}, Planetary and Space Science 78 (2013) 1 --
  21.
\newblock \href {http://dx.doi.org/https://doi.org/10.1016/j.pss.2012.12.002}
  {\path{doi:https://doi.org/10.1016/j.pss.2012.12.002}}.
\newline\urlprefix\url{http://www.sciencedirect.com/science/article/pii/S0032063312003777}

\bibitem{chebTaylor}
J.~D. Mireles~James, M.~Murray,
  \href{https://doi.org/10.1142/S0218127417300506}{Chebyshev-{T}aylor
  parameterization of stable/unstable manifolds for periodic orbits:
  Implementation and applications}, International Journal of Bifurcation and
  Chaos 27~(14) (2017) 1730050.
\newblock \href
  {http://arxiv.org/abs/https://doi.org/10.1142/S0218127417300506}
  {\path{arXiv:https://doi.org/10.1142/S0218127417300506}}, \href
  {http://dx.doi.org/10.1142/S0218127417300506}
  {\path{doi:10.1142/S0218127417300506}}.
\newline\urlprefix\url{https://doi.org/10.1142/S0218127417300506}

\bibitem{zhang}
L.~{Zhang}, R.~{de la Llave}, {Transition state theory with quasi-periodic
  forcing}, Communications in Nonlinear Science and Numerical Simulations 62
  (2018) 229--243.
\newblock \href {http://dx.doi.org/10.1016/j.cnsns.2018.02.014}
  {\path{doi:10.1016/j.cnsns.2018.02.014}}.

\end{thebibliography}
\bibliographystyle{elsarticle-num}

\end{document}